\documentclass{compositio}
\usepackage{amsopn,amsmath,amssymb,amscd,latexsym,textcomp}
\usepackage[all]{xy}
\newcommand{\A}{{\mathbb{A}}}

\newcommand{\C}{{\mathbb{C}}}
\newcommand{\Cc}{{\cal{C}}}
\newcommand{\D}{{\mathbb{D}}}
\newcommand{\F}{{\mathbb{F}}}

\newcommand{\N}{{\mathbb{N}}}
\newcommand{\Pe}{{\mathbb{P}}}
\newcommand{\Q}{{\mathbb{Q}}}
\newcommand{\oQ}{\overline{\Q}}
\newcommand{\opi}{\overline{\pi}}
\newcommand{\Rc}{{\cal R}}
\newcommand{\uT}{\underline{T}}

\newcommand{\Z}{{\mathbb{Z}}}
\newcommand{\ddet}{\mathrm{det}}

\newcommand{\dlog}{\mathrm{dlog}}
\newcommand{\diag}{\mathrm{diag}}
\newcommand{\ch}{\mathrm{ch}}
\newcommand{\eexp}{\mathrm{exp}}
\newcommand{\tr}{\mathrm{tr}}

\newcommand{\Gal}{\mathrm{Gal}}
\newcommand{\ggcd}{\mathrm{gcd}}

\newcommand{\ok}{\overline{k}}
\newcommand{\lcm}{\mathrm{lcm}}
\newcommand{\rk}{\mathrm{rk}}
\newcommand{\ddim}{\mathrm{dim}\;}
\newcommand{\Aut}{\mathrm{Aut}\,}

\newcommand{\ualpha}{\underline{\alpha}}
\newcommand{\End}{\mathrm{End}\,}

\newcommand{\Gl}{\mathrm{Gl}\,}
\newcommand{\Gr}{\mathrm{Gr}}
\newcommand{\Ind}{\mathrm{Ind}}
\newcommand{\kker}{\mathrm{ker}\,}
\newcommand{\Rep}{\mathrm{Rep}\,}
\newcommand{\uRep}{\underline{\mathrm{Rep}\,}}

\newcommand{\Res}{\mathrm{Res}}
\newcommand{\xa}{\mathrm{Res}_{G_L}^{G_k}}
\newcommand{\xb}{\mathrm{Ind}_{G_L}^{G_k}}

\newcommand{\za}{\mathrm{Rep}_{G_k} \Q_l}
\newcommand{\zb}{\mathrm{Rep}_{G_L} \Q_l}
\newcommand{\zc}{\mathrm{WRep}_{G_k} \Q_l}

\newcommand{\Spec}{\mathrm{Spec}}
\newcommand{\Var}{\mathrm{Var}}
\newcommand{\Vect}{\mathrm{Vect}}
\newcommand{\Oh}{{\mathcal O}}
\newcommand{\Rh}{{\mathcal R}}
\newcommand{\oX}{\overline{X}}
\newcommand{\oY}{\overline{Y}}
\newcommand{\dkr}{\mbox{\rm\textrbrackdbl}}
\newcommand{\dkl}{\mbox{\rm\textlbrackdbl}}
\newcommand{\dis}{\displaystyle}
\newtheorem{theorem}{Theorem}
\newtheorem{lemma}[theorem]{Lemma}
\newtheorem{prop}[theorem]{Proposition}
\newtheorem{defn}[theorem]{Definition}
\newtheorem{cor}[theorem]{Corollary}
\newenvironment{rem}{\noindent {\bf Remark}}{}

\newenvironment{proofof}{\noindent {\bf Proof of}}{\mbox{}\hfill$\Box$}

\begin{document}


\title[Algebraic independence]{Algebraic independence in the Grothendieck ring of varieties}
\author{N. Naumann}
\address{NWF I- Mathematik\\ Universit\"at Regensburg\\ Universit\"atsstrasse 31\\93053 Regensburg}

\begin{abstract}
We give sufficient cohomological criteria for the classes of
given varieties over a field $k$ to be algebraically independent in the 
Grothendieck ring of varieties over $k$ and construct some examples.
\end{abstract}

\maketitle

\section{Introduction}\label{introduction} 

Let $k$ be a field. The Grothendieck ring of varieties over $k$, denoted
$K_0(\Var_k)$, was defined by A. Grothendieck in \cite{Grothletter}:\\
The abelian group $K_0(\Var_k)$ is generated by the isomorphism classes
$[X]$ of separated $k$-schemes of finite type $X/k$ subject to the relations
$[X]=[Y]+[X-Y]$ for $Y\subset X$ a closed subscheme. Multiplication is 
given by $[X_1]\cdot [X_2]:=[X_1 \times_k X_2]$ and makes $K_0(\Var_k)$
a commutative ring with unit $[\Spec (k)]$.
By its very definition $K_0(\Var_k)$ is the value group of the universal
Euler-Poincar\'e characteristic with compact support for varieties over 
$k$ and is thus a fundamental invariant of the algebraic geometry over $k$.
The most outstanding result on the structure of $K_0(\Var_k)$ for 
$k$ of characteristic zero is a presentation of the ring $K_0(\Var_k)$
in terms of generators and relations anticipated by E. Looijenga and 
proved by his student F. Bittner \cite{Bittner}. B. Poonen showed 
that for $k$ of characteristic zero a linear combination of
the classes of suitable abelian varieties is a zero divisor in $K_0(\Var_k)$
\cite{Poonen} and J. Koll\'ar has computed the subring of $K_0(\Var_k)$
generated by the classes of conics for suitable base fields $k$ \cite{Kollar}.
A deep problem pertaining to the structure of $K_0(\Var_k)$
is the rationality of motivic zeta-functions posed by M. Kapranov \cite{Kapranov} on which there has been recent progress due to M. Larsen and
V. Lunts \cite{LaLu}. This problem is intimately related 
to the finite-dimensionality of motives, c.f. \cite{Andre} for an exposition.Furthermore, the ring $K_0(\Var_k)$ plays a central 
r\^ole in the theory of motivic integration, c.f. \cite{Looijengamotmeas}.
The present note is dedicated to the following problem about the structure 
of $K_0(\Var_k)$.\\
Given varieties $X_i/k$, when are their classes $[X_i]\in K_0(\Var_k)$
algebraically independent (in the sense made precise at the end of this
introduction) ?\\
We give a number of sufficient cohomological conditions for this to be the case and construct
some examples. Our main results are also valid in positive characteristic
where there is no known theorem about the structure of $K_0(\Var_k)$.
We now review the content of the individual sections.\\
In section \ref{a motivic measure} we construct various ring homomorphisms
$\mu_k:K_0(\Var_k)\longrightarrow R$ (so called motivic measures) all based
on \'etale cohomology. The weight filtration is used to show that the 
classes of proper smooth varieties over a finitely generated field $k$ are
algebraically independent in $K_0(\Var_k)$ as soon as their first \'etale
cohomology groups are algebraically independent in the ring of virtual 
$l$-adic Galois representations (corollary \ref{H1suffices}). We then reformulate
known properties of rationality and $l$-independence for varieties over 
$\F_q$ (proposition \ref{rationalimage}).\\
In subsection \ref{virtual continuous representations} we reduce the problem 
of algebraic independence of virtual $l$-adic Galois representations to a 
problem about representations of (possibly non-connected) reductive
groups, give a result illustrating the subtleties arising from this non-
connectedness (theorem \ref{abelian surfaces}) and construct an infinite
sequence of curves over a given finite field $k$ the classes of which are
algebraically independent in $K_0(\Var_k)$ (theorems \ref{inftycurvesfinitefield} and \ref{DHcurvesdo}). This shows in particular that
$K_0(\Var_k)$ contains a polynomial ring in infinitely many variables
and supports the intuition that $K_0(\Var_k)$ should be a large ring.
We obtain a similar result for number fields using only elliptic curves 
(theorem \ref{ellovernumber}).\\
In subsection \ref{using a result of skolem} we use a lemma of Skolem together
with the aforementioned rationality properties to give another approach 
to the algebraic independence of virtual $l$-adic Galois representations 
and show that a generic pair of curves over a finite field has algebraically independent classes in the Grothendieck ring of varieties
(theorem \ref{indepoftwocurves} in subsection \ref{twocurves}). We also determine the structure of the subring 
of $K_0(\Var_{\F_q})$ generated by zero dimensional varieties showing in
particular that $K_0(\Var_{\F_q})$ contains infinitely many zero divisors
(theorem \ref{structureofS}).

\begin{acknowledgements} I would like to thank my advisor C. Deninger 
for constant support and encouragement during the preparation 
of my thesis the main results of which are presented here. Furthermore,
I would like to thank H. Frommer, U. Jannsen, M. Kisin, B. Moonen, D. Roessler, C. Serp\'e
and C. Soul\'e for useful discussions and the referee for a very thorough 
report.
\end{acknowledgements}

The following definition will be central throughout the rest of this note.
\begin{defn} Let $R$ be a commutative unitary ring and let $n\Z$ be the 
kernel of $\Z\longrightarrow R$. We call a set $\{x_i: i\in I\}$ of elements
$x_i\in R$ algebraically independent if the morphism

\[
(\Z/n)[T_i:i\in I]\longrightarrow R\; , \; T_i\mapsto x_i
\]

is injective.
\end{defn}

If $R$ is a field then this means algebraically independent over the prime
filed of $R$ in the usual sense. Note that $1/2\in\Q$ is algebraically
dependent though not integral over $\Z$, c.f. theorem \ref{zerodimthm}.
All the rings appearing below will have $n=0$, i.e. $\Z\subset R$. In this
case, elements $x_i\in R$ are algebraically independent if and only if 
the subring of $R$ they generate is a polynomial ring over $\Z$ in the
variables $x_i$.\\

\section{A motivic measure}\label{a motivic measure}

\subsection{General construction} \label{general construction}

Let $k$ be a field. We fix a separable closure $\ok$ of $k$ and a rational
prime $l$ different from the characteristic of $k$. For $X/k$ separated 
and of finite type we denote by

\[ 
H_c^i(\oX):=H_c^i(X\times_k\ok,\Q_l)\mbox{ for }i\ge 0
\]

the \'etale cohomology with compact support and constant coefficients $\Q_l$ of 
the base change $X\times_k\ok$. Then $H_c^i(\oX)$ is naturally a $G_k:=
\Aut(\ok/k)$-module and denoting by $\Rep_{G_k}\Q_l$ the category of 
finite dimensional continuous representations of $G_k$ over $\Q_l$
(i.e. $l$-adic Galois representations) we have a motivic measure

\begin{equation}\label{one}
\mu_k:K_0(\Var_k)\longrightarrow K_0(\Rep_{G_k}\Q_l)\; , \; [X]\mapsto\sum_i(-1)^i[H_c^i(\oX)].
\end{equation}

In fact, (\ref{one}) is (well-defined and) a homomorphism of abelian groups 
by excision, multiplicative by the K\"unneth formula and clearly preserves
the unit.\\
For a finite extension $k\subset L\subset\ok$ we have a ring homomorphism
(base change)
\[
-\times_k L: K_0(\Var_k)\longrightarrow K_0(\Var_L)\; , \; [X]\mapsto [X\times_k L]
\]
and a homomorphism of abelian groups (restriction of scalars)
\[
/k: K_0(\Var_L)\longrightarrow K_0(\Var_k)\; , \; [X\longrightarrow\Spec(L)]\mapsto [X\longrightarrow\Spec(L)\longrightarrow\Spec(k)].
\]
From the inclusion $G_L\subset G_k$ we have an exact tensor functor (restriction)
\begin{equation}\label{alpha}
\Res_{G_L}^{G_k}:\Rep_{G_k}\Q_l\longrightarrow\Rep_{G_L}\Q_l
\end{equation}
and an exact functor (induction)
\begin{equation}\label{beta}
\Ind_{G_L}^{G_k}:\Rep_{G_L}\Q_l\longrightarrow\Rep_{G_k}\Q_l\; , \; V\mapsto\Q_l[G_k]\otimes_{\Q_l[G_L]} V
\end{equation}
inducing a ring homomorphism and a homomorphism of abelian groups on the level of $K_0$ to be denoted
by the same symbols. The following diagram is commutative:
\[
\xymatrix{ K_0(\Var_k) \ar_{-\times_k L}[d] \ar^-{\mu_k}[r] & K_0(\za) \ar_{\xa}[d] \\
K_0(\Var_L) \ar@<-1ex>_{/k}[u] \ar^-{\mu_L}[r] & K_0(\zb). \ar@<-1ex>_{\xb}[u] }
\]

\subsection{Weight filtration}\label{weight filtration}

Notation being as in subsection \ref{general construction} we now assume in
addition that $k$ is finitely generated. In this subsection we will
incorporate the weight filtration on $H_c^*(\oX)$ into the motivic
measure (\ref{one}). We need to assume that $k$ is finitely generated 
in order to have a theory of weights. Vaguely speaking, this is a substitute
(valid in any characteristic) of Hodge theory and will allow us to isolate
$H^1$ from the Euler-Poincar\'e characteristic in theorem \ref{H1} below.\\
We refer the reader to \cite{Jannsen},\S 6 for the definition of the full
subcategory

\[
\zc\subset\za
\]

of $l$-adic Galois representations having a weight filtration. From \cite{Jannsen},
lemma 6.8.2 and equation (6.8.3) we know that $H_c^i(\oX)\in\zc$ holds
for all $X/k$ separated and of finite type and hence (\ref{one}) factors as
\begin{equation}\label{fourandhalf}
\mu_k:K_0(\Var_k)\longrightarrow K_0(\zc).
\end{equation}
For $V\in\zc$ we denote by $W_{\bullet}V$ the weight filtration of $V$
and by $\Gr_i^W(V):=W_iV/W_{i-1}V$ $(i\in\Z)$ the associated graded of weight
$i$. As any morphism in $\zc$ is strictly compatible with the weight filtration 
(\cite{Jannsen}, lemma 6.8.1) the functor

\[ 
\Gr_i^W:\zc\longrightarrow\za\; , \; V\mapsto\Gr_i^W(V)
\]

is exact (\cite{DeligneHodge2}, proposition 1.1.11, ii)) and thus induces 
a homomorphism 

\[ 
\Gr_i^W:K_0(\zc)\longrightarrow K_0(\za).
\]

So we have a homomorphism of abelian groups 

\begin{equation}\label{four}
\Phi: K_0(\zc)\longrightarrow K_0(\za)[T,T^{-1}]\; , \; x\mapsto\sum_{i\in\Z}\Gr_i^W(x)T^i
\end{equation}

which is a ring homomorphism as follows from

\[
\Gr_i^W(V\otimes W)\simeq\oplus_{a+b=i}\Gr_a^W(V)\otimes\Gr_b^W(W)\mbox{ for }
i\in\Z\mbox{ and }V,W\in\zc.
\]

As $\Gr_i^W(H^*_c(\oX))=0$ for $i<0$ (see \cite{Katzreview}) we obtain, composing (\ref{fourandhalf}) and (\ref{four}), a motivic measure

\begin{equation}\label{five}
\mu_k:K_0(\Var_k)\longrightarrow K_0(\za)[T]
\end{equation}

given explicitly by

\begin{equation}\label{six}
\mu_k([X])=\sum_{i\ge0}\left( \sum_j (-1)^j[\Gr_i^W(H^j_c(\oX))]\right)T^i.
\end{equation}

The slight abuse of notation in denoting (\ref{one}) and (\ref{five})
by the same symbol $\mu_k$ will cause no confusion.\\
For a finite extension $k\subset L\subset\ok$ we extend
the morphisms induced on $K_0$ by (\ref{alpha})
(resp. (\ref{beta})) to $K_0(\za)[T]$ (resp. $K_0(\zb)[T]$) by demanding that
$T\mapsto T$. With $\mu_k$ as in (\ref{five}) we then have a commutative diagram

\begin{equation}\label{seven}
\xymatrix{ K_0(\Var_k) \ar_{-\times_k L}[d] \ar^-{\mu_k}[r] & K_0(\za)[T] \ar_{\xa}[d] \\
K_0(\Var_L) \ar@<-1ex>_{/k}[u] \ar^-{\mu_L}[r] & K_0(\zb)[T]. \ar@<-1ex>_{\xb}[u] }
\end{equation}

The commutativity of (\ref{seven}) follows from the compatibility of the weight
filtration with restriction and induction which we leave to the reader to verify
using the uniqueness of the weight filtration (\cite{Jannsen}, lemma 6.8.1, a)).
We now use the weight filtration to establish the following criterion for 
algebraic independence in $K_0(\Var_k)$:\\
For $X / k$ separated and of finite type and $i\ge 0$ we write
\[
\Gr_i^W (X) := \sum_j (-1)^{j} [\Gr_i^W (H^j_c (\oX))] \in K_0 (\Rep_{G_k} \Q_l) \; .
\]
Then $\mu_k ([X]) = \sum_i \Gr_i^W (X) T^i\in K_0(\Rep_{G_k}\Q_l)[T]$, see (\ref{six}).

\begin{theorem}
  \label{H1}
Let $k$ be finitely generated and let $X_1 , \ldots , X_n / k$ be separated and of finite type and assume that $\Gr^W_0 (X_i) \in \Z \subseteq K_0 (\Rep_{G_k} \Q_l)$ for all $i$ and that the $\Gr^W_1 (X_i)$ are algebraically independent in $K_0 (\Rep_{G_k} \Q_l)$. Then the $\mu_k ([X_i])$ are algebraically independent in $K_0 (\Rep_{G_k} \Q_l) [T]$ and hence the $[X_i]$ are algebraically independent in $K_0(\Var_k)$. 
\end{theorem}
Note that $K_0 (\Rep_{G_k} \Q_l)$ is augmented over $\Z$ by the degree map, hence $\Z \subseteq K_0 (\Rep_{G_k} \Q_l)$. The last assertion of 
the theorem follows from the existence of the motivic measure (\ref{five}).

\begin{cor}\label{H1suffices}
  Let $k$ be finitely generated and let $X_1 , \ldots , X_n / k$ be proper and smooth. Then, if the $[H^1_c (\oX_i)]$ are algebraically independent in $K_0(\Rep_{G_k}\Q_l)$, so are the $\mu_k ([X_i])$ in $K_0(\Rep_{G_k}\Q_l)[T]$ and hence also the $[X_i]$ in
$K_0(\Var_k)$. 
\end{cor}

\begin{proof} 
As the $X_i/k$ are proper and smooth we have $\Gr^W_0(X_i)=[H_c^0(\oX_i)]
=[\Q_l^{\oplus |\pi_0(\oX_i)|}]=|\pi_0(\oX_i)|\in\Z$ and $\Gr_1^W(X_i)=-[H^1_c(\oX_i)]$. Since the $[H^1_c(\oX_i)]$ are algebraically independent by assumption, so 
are the $\Gr_1^W(X_i)$ and theorem \ref{H1} applies.
\end{proof}

For the proof of theorem \ref{H1} we need some elementary preparation.
Let $R$ be a commutative unitary ring. On $R [T_1 , \ldots , T_n]$ we have the usual derivations $\frac{\partial}{\partial T_i}$ and for $\ualpha = (\alpha_1 , \ldots , \alpha_n) \in \N^n_0$ we write $\frac{\partial^{|\ualpha|}}{\partial T^{\alpha_1}_1 \ldots \partial T^{\alpha_n}_n}$ for their iterations; here $|\ualpha| := \sum \alpha_i$.\\
Let now be given $f_1 , \ldots , f_n \in R [T] , G \in \Z [T_1 , \ldots , T_n]$ and put $\tilde{G} (T) := G (f_1 , \ldots , f_n) \in R [T]$. For $N \ge 1$ let $I_N \subseteq R [T]$ be the ideal generated by $\frac{\partial^{|\ualpha|} G}{\partial T^{\alpha_1}_1 \ldots \partial T^{\alpha_n}_n} (f_1 , \ldots , f_n)$ for $|\ualpha| \le N$ and put $I_0 := 0$. The $(I_N)_{N\geq 0}$ form an ascending
chain of ideals of $R[T]$ and $\frac{d}{dT}(I_N)\subset I_{N+1}$ for $N\ge 0$.

\begin{prop}
  \label{calculus}
For $N \ge 1$ we have
\[
\frac{d^N \tilde{G}}{d T^N} \equiv \sum_{1 \le  i_1 , \ldots , i_N \le n} \frac{\partial^N G}{\partial T_{i_1} \ldots \partial T_{i_N}} (f_1 , \ldots , f_n) \frac{d f_{i_1}}{dT} \ldots \frac{df_{i_N}}{dT} \; \mbox{ (mod } I_{N-1})\; .
\]
\end{prop}

\begin{proof}
  We use induction on $N \ge 1$. We have
\[
\frac{d\tilde{G}}{dT} = \frac{d}{dT} G (f_1 , \ldots , f_n) = \sum_{1 \le i_1 \le n} \frac{\partial G}{\partial T_{i_1}} (f_1 , \ldots , f_n) \frac{df_{i_1}}{dT} \; ,
\]
i.e. the result for $N = 1$.
We assume $N \ge 2$ and compute, using the induction hypothesis:
\[
\frac{d^N \tilde{G}}{dT^N} = \frac{d}{dT} \left( \sum_{1 \le i_1 , \ldots , i_{N-1} \le n} \frac{\partial^{N-1} G}{\partial T_{i_1} \ldots \partial T_{i_{N-1}}} (f_1 , \ldots , f_n) \frac{df_{i_1}}{dT} \cdots \frac{df_{i_{N-1}}}{dT} + \alpha \right) \; ,
\]
some $\alpha \in I_{N-2}$. As $\frac{d}{dT} (I_{N-2}) \subseteq I_{N-1}$ this equals, modulo $I_{N-1}$:
\begin{eqnarray*}
  \lefteqn{\sum_{1 \le i_1 , \ldots , i_{N-1} \le n} \left[ \left( \sum^n_{j=1} \frac{\partial^N G}{\partial T_{i_1} \ldots \partial T_{i_{N-1}} \partial T_j} (f_1 , \ldots , f_n) \frac{d f_j}{dT} \right) \frac{df_{i_1}}{dT} \cdots \frac{df_{i_{N-1}}}{dT} \right. }\\
 & + & \left. \frac{\partial^{N-1} G}{\partial T_{i_1} \ldots \partial T_{i_{N-1}}} (f_1 , \ldots , f_n) \frac{d}{dT} (...) \right] \overset{i_N := j}{\equiv} \\
\lefteqn{\sum_{1 \le i_1 , \ldots , i_N \le n} \frac{\partial^N G}{\partial T_{i_1} \ldots \partial T_{i_N}} (f_1 , \ldots , f_n) \frac{df_{i_1}}{dT} \cdots \frac{df_{i_N}}{dT} \; , }
\end{eqnarray*}

as claimed. 

\end{proof}

\begin{proofof} theorem \ref{H1}. Put $R := K_0 (\Rep_{G_k} \Q_l)$ ,\linebreak $f_i(T) := \mu_k ([X_i]) = \sum_{\nu \ge 0} \Gr_{\nu}^W (X_i) T^{\nu} \in R [T] , a_i := \Gr_0^W (X_i) \in \Z\subset R$ by assumption. Assume that we have $G\in \Z [T_1 , \ldots , T_n]$ with $G (\mu_k ([X_1]), \ldots, \mu_k([X_n])) = 0$, i.e. 
  \begin{equation}
    \label{built}
    \tilde{G} (T) := G (f_1 (T) , \ldots , f_n (T)) = 0 \; .
  \end{equation}
We have to show that $G = 0$ which follows from the\\
{\it Claim}: For $N \ge 0$ and any $\ualpha \in \N^n_0$ with $|\ualpha| \le N$ we have
\[
\frac{\partial^{|\ualpha|}G}{\partial T^{\alpha_1}_1 \ldots \partial T^{\alpha_n}_n} (a_1 , \ldots , a_n) = 0 \; ,
\]
which in turn will be established by induction on $N$.
Note that $f_i (0) = a_i$, so from (\ref{built}) putting $T = 0$ we get
$G (a_1 , \ldots , a_n) = 0$,
i.e. the claim for $N = 0$. Now assume that $N \ge 1$. By applying proposition \ref{calculus} to (\ref{built}) and putting $T = 0$ we obtain the following relation in $R$:
\[
0 = \frac{d^N \tilde{G}}{dT^N} (T = 0) = \sum_{1 \le i_1 , \ldots , i_N \le n} \frac{\partial^N G}{\partial T_{i_1} \ldots \partial T_{i_N}} (a_1 , \ldots , a_n) \frac{df_{i_1}}{dT} (0) \ldots \frac{df_{i_N}}{dT} (0) + \alpha (0) \; .
\]
Here $\alpha \in I_{N-1}$, hence $\alpha (0) = 0$ by induction hypothesis.
As $\frac{df_i}{dT} (0) = \Gr_1^W (X_i)$ we have the following relation in $R$:
\begin{equation}
  \label{troll}
  0 = \sum_{1 \le i_1 , \ldots , i_N \le n} \frac{\partial^N G}{\partial T_{i_1} \ldots \partial T_{i_N}} (a_1 , \ldots , a_n) \Gr_1^W (X_{i_1}) \ldots \Gr_1^W (X_{i_N}) \; .
\end{equation}
We collect terms in this expression: Consider
\begin{eqnarray*}
  \pi : I := \{ (i_1 , \ldots , i_N) : 1 \le i_j \le n \} & \longrightarrow & \{ \ualpha=(\alpha_1 , \ldots , \alpha_n)\in \N_0^n : |\ualpha|=N \} \\
(i_1 , \ldots , i_N) & \longmapsto & (\alpha_k := | \{ 1 \le j \le N : i_j = k \} | )_{k = 1 , \ldots , n} \; .
\end{eqnarray*}
Given $(i_1 , \ldots , i_N)\in I$ with $\pi ((i_1 , \ldots , i_N)) = (\alpha_1 , \ldots , \alpha_n)$ we have
\[
\begin{array}{l}
\dis \frac{\partial^N G}{\partial T_{i_{1}} \ldots \partial T_{i_N}} = \frac{\partial^N G}{\partial T^{\alpha_1}_1 \ldots \partial T^{\alpha_n}_n} \quad \mbox{and} \vspace{0.5cm} \\

\dis \Gr_1^W (X_{i_1}) \ldots \Gr_1^W (X_{i_N}) = \Gr_1^W (X_1)^{\alpha_1} \ldots \Gr_1^W (X_n)^{\alpha_n} \; .
\end{array}
\]
So (\ref{troll}) may be written as
\begin{equation}
  \label{troi}
  \sum_{\ualpha = (\alpha_1 \ldots \alpha_n) \in \N_0^n \atop |\ualpha| = N} |\pi^{-1} (\ualpha)| \frac{\partial^N G}{\partial T^{\alpha_1}_1 \ldots \partial T^{\alpha_n}_n} (a_1 , \ldots , a_n) \Gr_1^W (X_1)^{\alpha_1} \ldots \Gr_1^W (X_n)^{\alpha_n} = 0 \; .
\end{equation}
The $\Gr_1^W (X_i) \in R$ are algebraically independent by assumption and $\pi$ is surjective, hence $|\pi^{-1} (\ualpha)| \neq 0$. So (\ref{troi}) implies
\[
\frac{\partial^N G}{\partial T^{\alpha_1}_1 \ldots \partial T^{\alpha_n}_n} (a_1 , \ldots , a_n) = 0
\]
for all $\ualpha \in \N_0^n$ with $|\ualpha| = N$, concluding the proof of the claim and of theorem \ref{H1}. 
\end{proofof}

\subsection{Finite base field}\label{finite base field}

In the situation of \ref{general construction} we now assume in addition that
the base field $k=\F_q$ is a finite field.
We will use the fact that $G_{\F_q}$ is topologically cyclic, generated by 
the geometric Frobenius $F_q\in G_{\F_q}$, to rewrite the results of \ref{general construction}
and \ref{weight filtration} in terms of characteristic polynomials of 
Frobenius.
We need to recall some facts about the ``universal ring'' from \cite{DG},
V, \S 5, no 2:\\
For any commutative unitary ring $R$ the set $\Lambda(R):=1+tR\dkl t \dkr$
of formal power series with coefficients in $R$ and constant coefficient 1 is
an abelian group under multiplication. The resulting functor $R\mapsto\Lambda(R)$ can be endowed with the structure of a functor in commutative
unitary rings such that for all $R$ and $a,b\in R$

\begin{equation}\label{1}
(1-at)(1-bt)=1-(ab)t\mbox{ in }\Lambda(R)=1+t R\dkl t\dkr,
\end{equation}

see {\em loc. cit.} 2.2. The multiplicative unit of $\Lambda(R)$ is $1-t$.
The ring structure of $\Lambda(R)$ may cause confusion:
The addition is the usual multiplication of power series but the multiplication
in $\Lambda(R)$ is a rather unusual operation on power series. We will
denote addition and multiplication in $\Lambda(R)$ by the usual symbols
but make clear that the composition is to be understood in $\Lambda(R)$.
In this sense, the additive analogue of (\ref{1}) is 
$(1-at)+(1-bt)=1-(a+b)t+abt^2$ in $\Lambda(R)$.\\
Let $R$ and $f\in\Lambda(R)$ be given. One defines $\partial_{\nu}(f)\in R$
for all $\nu\ge 1$ by the expansion $\dlog (f)=:-\sum_{\nu\ge 1}\partial_{\nu}(f)t^{\nu -1}$. Here, $\dlog(f)=\frac{df/dt}{f}$ (usual quotient of power series) is the 
logarithmic derivative of $f$. If we denote by $\Oh$ the identity functor
on the category of commutative unitary rings then each $\partial_{\nu}:\Lambda\longrightarrow\Oh$ is a natural transformation ({\em loc. cit.} 2.4) and if $R$ is a $\Q$-algebra

\begin{equation}\label{splitoverQ}
\Lambda(R)\stackrel{\sim}{\longrightarrow}R^{\N}\; , \; f\mapsto
(\partial_{\nu}(f))_{\nu\ge 1}
\end{equation}

is an isomorphism of rings in which $R^{\N}$ is the usual product of countably
many copies of the ring $R$.
For any $\nu\ge 1$ there is a natural transformation of ring valued functors
(``Frobenius'') $F_{\nu}:\Lambda\longrightarrow\Lambda$ characterized by

\begin{equation}\label{3}
F_{\nu}(1-at)=1-a^{\nu}t\mbox{ in }\Lambda(R)=1+t R\dkl t \dkr
\end{equation}

for all $R$ and $a\in R$ and a natural transformation of functors in abelian 
groups (``Verschiebung'') $V_{\nu}:\Lambda\longrightarrow\Lambda$ characterized by $V_{\nu}(1-at)=1-at^{\nu}\mbox{ in }\Lambda(R)$.
In fact, we have

\begin{equation}\label{4}
V_{\nu}(f(t))=f(t^{\nu})\mbox{ for all }f(t)\in\Lambda(R)=1+tR\dkl t \dkr,
\end{equation}

see {\em loc. cit.} 2.5.\\
The rings $\Lambda(R)$ are convenient tools to keep track of 
characteristic polynomials:

\begin{prop}\label{charcount}
Let $k$ be a field and $V,W$ finite dimensional $k$-vectorspaces.

\begin{tabular}{cl}
i) & For $\phi\in\End_k(V)$ and $\psi\in\End_k(W)$ we have\\
 & \\
   & $\ddet(1-(\phi\oplus\psi)t|V\oplus W)=\ddet(1-\phi t|V)+\ddet(1-\psi t|W)\mbox{ and}$\\
 & $\ddet(1-(\phi\otimes\psi)t|V\otimes W)=\ddet(1-\phi t|V) \ddet(1-\psi t|W)\mbox{ in }\Lambda(k)=1+tk\dkl t \dkr.$\\
 & \\
 & For any $\nu\ge 1$:\\
 & \\
 & $\ddet(1-\phi^{\nu} t|V)=F_{\nu}(\ddet(1-\phi t|V))\mbox{ in }\Lambda(k)$ and \\
 & $\partial_{\nu}(\ddet(1-\phi t|V))=\tr(\phi^{\nu}|V)$ in $k$.\\
 & \\

ii) & Let $G$ be a group and $N\subset G$ a normal subgroup such that $G/N$\\
 & is finite cyclic of order $\nu$ generated by the class of $\pi\in G$.\\
 & If $V$ is a $k[N]$-module, finite dimensional as a $k$-vectorspace, \\
 & then we have
\end{tabular}

\begin{equation}\label{2}
\ddet(1-\pi t|\Ind_N^G(V))=V_{\nu}(\ddet(1-\pi^{\nu} t|V))\mbox{ in }\Lambda(k)=1+t k\dkl t \dkr.
\end{equation}

\end{prop}

In $ii)$ we denoted by $\Ind_N^G(V):=k[G]\otimes_{k[N]} V$ the $k[G]$-
module induced by $V$. Note that $\pi^{\nu}\in N$, hence the right hand side
of (\ref{2}) is meaningful.

\begin{proof}: Part i) follows from equations (\ref{1})
and (\ref{3}) and \cite{Hartshorne}, appendix C, lemma 4.1. Using (\ref{4}), part ii)
is equivalent to 
$\ddet(1-\pi t|\Ind_N^G(V))=\ddet(1-\pi^{\nu} t^{\nu}|V)$.
This has to be checked by an explicit computation which we omit.

\end{proof}

For a ring $R$ we define a subset $\Rh(R)\subset\Lambda(R)$ by

\[
\Rh(R):=\{ f-g\;|\; f,g\in 1+tR[t]\}\subset\Lambda(R)=1+t R\dkl t \dkr.
\]

It is clear that $\Rh\subset\Lambda$ is a subfunctor in abelian groups.
One should think of the elements of $\Rh(R)$ as the the power series expansions
of certain rational functions, namely the quotients of polynomials with 
constant coefficient 1. 

\begin{prop}\label{defcurlyr}
$\Rh\subset\Lambda$ is a subfunctor in rings.
\end{prop}

\begin{proof} 
Firstly, $1_{\Lambda(R)}=1-t\in\Rh(R)\subset\Lambda(R)=1+t R \dkl t \dkr$. Now it suffices to check that for 
$f,g\in 1+tR[t]$ we have $fg\in\Rh(R)\subset\Lambda(R)$. Indeed, $fg$ is again a polynomial according to \cite{DG}, p. 632, formula (**).
\end{proof}

Using proposition \ref{charcount} one sees that there is a unique ring homomorphism

\[
\chi_q: K_0(\Rep_{G_{\F_q}}\Q_l)\hookrightarrow \Lambda(\Q_l)
\]

such that

\[
\chi_q([V])=\det(1-F_q t|V)\in\Lambda(\Q_l)=1+t\Q_l\dkl t \dkr\mbox{ for } V\in\Rep_{G_{\F_q}}\Q_l
\]

and $\chi_q$ is injective as follows from \cite{Bourbakialgebra}, \S 12,1, proposition 3.\\
Composing (\ref{five}) with $\chi_q$ we obtain a motivic measure

\begin{equation}\label{30}
\mu_q:K_0(\Var_{\F_q})\longrightarrow\Lambda(\Q_l)[T]
\end{equation}

given explicitly by

\begin{equation}\label{formula5}
\mu_q([X])=\sum_{i\ge 0}\left(\sum_j(-1)^j\ddet(1-F_q t |\Gr_i^W(H^j_c(\oX)))\right)T^i.
\end{equation}

\begin{prop}\label{rationalimage}
The image of (\ref{30}) is contained in 
$\Rh(\Q)[T]\subset\Lambda(\Q_l)[T]$.
\end{prop}

\begin{proof} Let $X/\F_q$ be separated and of 
finite type. The cohomological formula for the zeta function $Z_{X/{\F_q}}(t)$
of $X/{\F_q}$ (\cite{Mi2}, theorem 13.1) gives in $\Lambda(\Q_l)=1+t\Q_l\dkl t \dkr$:

\begin{eqnarray*}
Z_{X/{\F_q}}(t)  =  \sum_{j\ge 0}(-1)^{j+1}\ddet(1-F_q t|H^j_c(\oX))= & & \\
\sum_{j\ge 0}(-1)^{j+1}\left(\sum_i\ddet(1-F_q t|\Gr_i^W(H^j_c(\oX)))\right) = &  & \\
\sum_i\left(\sum_{j\ge 0}(-1)^{j+1}\ddet(1-F_q t|\Gr_i^W(H^j_c(\oX)))\right).& &  
\end{eqnarray*}

In this expression, every summand is in $\Q_l(t)$ and we have that $Z_{X/\F_q}(t)\in\Q\dkl t \dkr$ from the very definition of $Z_{X/\F_q}(t)$.
For weight reasons there can be no cancellation among different summands (recall that the sums are actually products of power series), hence every
summand lies in $\Q_l(t)\cap\Q\dkl t \dkr\subset\Q(t)$. As every summand
is in addition a quotient of polynomials with constant coefficient 1, every summand
lies in fact in $\Rh(\Q)$ and comparison with (\ref{formula5}) concludes the proof.
\end{proof}

The relevance of having a motivic measure taking values in $\Rh(\Q)[T]$ rather
than in $\Lambda(\Q_l)[T]$ will become clear in \ref{using a result of skolem}.
Finally we record some formulae relating the various operations on $\Lambda$
with the motivic measure $\mu_q$.

\begin{prop}\label{motmeasurefiniteproperties}

\[
\begin{array}{ll}
i) & \mbox{If } X/\F_q \mbox{ is separated and of finite type then, for all }\nu\ge 1,\\
 & \partial_{\nu}(\mu_q([X]))=\sum_{i\ge 0}\left(\sum_j(-1)^j\tr(F_q^{\nu}|\Gr_i^W(H^j_c(\oX)))\right)T^i\mbox{ in }\Q[T]. \\
 & \\
ii) & \mbox{For } \F_q\subset \F_{q^{\nu}}\subset\overline{\F_q} \mbox{ the finite extension of degree } \nu\ge 1 \\
 & \mbox{we have a commutative diagram}
\end{array}
\]

\begin{equation}\label{33}
\xymatrix{ K_0(\Var_{\F_q}) \ar_{-\times_{\F_q} \F_{q^{\nu}}}[d] \ar^-{\mu_q}[r] & \Rh(\Q)[T] \ar_{F_{\nu}}[d] \\
K_0(\Var_{\F_{q^{\nu}}}) \ar@<-1ex>_{/\F_q}[u] \ar^-{\mu_{q^{\nu}}}[r] & \Rh(\Q)[T]. \ar@<-1ex>_{V_{\nu}}[u] }
\end{equation}

\end{prop}

In $i)$ $\partial_{\nu}$ denotes the restriction of $\partial_{\nu}:\Lambda(\Q)\longrightarrow\Q$ to $\Rh(\Q)\subset\Lambda(\Q)$ extended to $\partial_{\nu}:\Rh(\Q)[T]\longrightarrow\Q[T]$ by demanding that $\partial_{\nu}(T)=T$ and in $ii)$ $F_{\nu}$
and $V_{\nu}$ denote the maps which on $\Rh(\Q)$ are the restrictions of 
$F_{\nu},V_{\nu}:\Lambda(\Q)\longrightarrow\Lambda(\Q)$ and with $F_{\nu}(T)=T$ and $V_{\nu}(T)=T$.

\begin{proof} Part $i)$ follows from (\ref{formula5}) and proposition \ref{charcount},i). Part $ii)$ follows from (\ref{seven}) and
proposition \ref{charcount},i) and ii).
\end{proof}

\section{Algebraic Independence}

\subsection{Virtual continuous representations}\label{virtual continuous representations}

In view of corollary \ref{H1suffices} we are led to consider the following problem:\\
Given a topological group $T$ and a finite extension $K/\Q_l$ denote by
$\Rep_T K$ the category of finite dimensional continuous representations of $T$ over $K$. Given finitely many $V_i\in\Rep_T K$, when are their classes 
$[V_i]\in K_0(\Rep_T K)$ algebraically independent ?\\
In this subsection we will use a Tannakian argument to reduce this question 
to a problem about algebraic independence of representations of a
(possibly non-connected) reductive algebraic group and in a special case
of a diagonalizable group. For this latter case we establish a rather explicit
Jacobi criterion (lemma \ref{jacobigross}). The subtleties arising from the possible non-connectedness of the reductive group
above will be illustrated by an example (theorem \ref{abelian surfaces}) and
motivate the following definition.

\begin{defn} Representations $V_i\in\Rep_T K$ are {\em geometrically} algebraically
independent if for all open subgroups of finite index $T'\subset T$ the 
classes of the restrictions $[\Res^T_{T'}(V_i)]\in K_0(\Rep_{T'}K)$ are algebraically independent.
\end{defn}

We now give a series of reduction steps to obtain workable criteria for the 
$[V_i]\in K_0(\Rep_T K)$ to be algebraically (in)dependent. The reader is
invited to glance at a simple application, theorem \ref{abelian surfaces}, first.\\
Denoting by $V_i^{ss}$ the semi-simplification of $V_i$ we have
\begin{equation}\label{eins}
[V_i]=[V_i^{ss}]\mbox{ in } K_0( \Rep_T K)
\end{equation}

and we can assume the $V_i$ to be semi-simple themselves.
In the sequel we will repeatedly use the obvious fact that elements
$x_i\in R\subset S$ of a subring $R$ of $S$ are algebraically independent
in $R$ if and only if they are so in $S$.
We now use a Tannakian argument to find an accessible subring
of $K_0( \Rep_T K)$ containing all the $[V_i]$. The reference for Tannakian categories is \cite{DeligneTannaka}.
With the obvious notion of tensor product, dual and using the forgetful
functor $\Rep_T K\longrightarrow\Vect_K$ to the category of 
finite dimensional $K$-vectorspaces as a fibre functor, the category  $\Rep_T K$
is a neutral Tannakian category over $K$.
We put $V:=\oplus_{i=1}^n V_i$ and consider

\[
i:\Cc:=< V >^{\otimes}\hookrightarrow \Rep_T K,
\]

the Tannakian subcategory generated by $V$, i.e. the smallest
strictly full subcategory of $\Rep_T K$ containing $V$ and stable under
the formation of finite direct sums, tensor products, duals, subobjects and quotients.
As $i$ is an exact tensor functor it induces a ring homomorphism

\begin{equation}\label{zwei}
K_0(\Cc)\hookrightarrow K_0( \Rep_T K)
\end{equation}

which is injective: one checks that $i$ preserves simple objects and
uses that the $K_0$ of an abelian category with all objects of finite length is $\Z$-free on the isomorphism classes of its simple objects \cite{Quillen}, \S 5, corollary 1. By construction we have $[V_i]\in K_0(\Cc)$
and our initial problem is reduced to considering the algebraic (in)dependence
of the $[V_i]$ as elements of $K_0(\Cc)$.
Let 
\begin{equation}\label{sternchen}
\phi:T\longrightarrow \Gl(V)(K)
\end{equation}

be the continuous homomorphism giving the action of $T$ on $V$. Here
$\Gl(V)$ denotes the algebraic group and $\Gl(V)(K)$ its group
of $K$-points endowed with the topology inherited from the one of $K$.
The Zariski closure of the image of $\phi$, say
\begin{equation}\label{drei}
G:=\overline{\mathrm{im}(\phi)}\subset\Gl(V),
\end{equation}
is an affine algebraic group over $K$ (\cite{Borel}, proposition 1.3, b)) and
can be identified with the Tannaka group of $\Cc$ as follows (\cite{SerreHodge-Tate}, 1.3).
From $\phi$ we have a continuous homomorphism $\tilde{\phi}:T\longrightarrow G(K)$
and the restriction of the fibre functor $\Cc\longrightarrow\Vect_K$
induces an equivalence of Tannakian categories

\begin{equation}\label{qwe}
\Cc\stackrel{\sim}{\longrightarrow}\uRep_K G
\end{equation}

where $\uRep_K G$ denotes the category of finite dimensional representations of the algebraic group $G$ over $K$. A quasi-inverse of (\ref{qwe}) is
given by

\[
\uRep_K G\stackrel{\sim}{\longrightarrow}\Cc\; , \; (G\stackrel{\psi}{\rightarrow}\Gl(W))\mapsto
\psi(K)\circ\tilde{\phi}.
\]

As $\Cc$ is semi-simple by (\ref{eins}) \cite{DMOS}, remark 2.28 implies
that $G$ is (not necessarily connected) reductive. The isomorphism $K_0(\Cc)\simeq K_0(\uRep_K G)$ induced by (\ref{qwe})
reduces our initial problem to studying the algebraic (in)dependence
of given representations of a reductive group.
We next show that we may allow finite extensions of the coefficients
$K$ of these representations.
For $K\subset L$ a finite extension, the exact tensor functor

\[
\uRep_K G\longrightarrow\uRep_L G\otimes L\; , \; V\mapsto V\otimes_K L
\]

induces a ring homomorphism

\begin{equation}\label{vier}
K_0(\uRep_K G)\hookrightarrow K_0(\uRep_L G\otimes L)
\end{equation}

which is injective, c.f. \cite{Milnemotivesoverfinite} lemma 3.14 for a similar argument. This allows us to replace $G$ by $G\otimes L$
and to assume that the connected component of identity $G^0\subset G$
is a connected split reductive group and that the finite \'etale group
$G/G^0$ is in fact constant.
Extending $K$ further to contain a splitting field of $G/G^0$ the ring
$K_0(\uRep_K G/G^0)$ is a familiar object from the representation theory of
finite (abstract)
groups \cite{Serrefinitegroups}, chapter II. Also $K_0(\uRep_K G^0)$
can be described quite precisely (c.f. below) but without imposing further, rather
restrictive, conditions (e.g. $G\simeq G^0\times G/G^0$) I cannot say much
about the structure of $K_0(\uRep_K G)$.
We will describe two special situations in which the above line of 
thought may be continued:\\
So $K/\Q_l$ is a finite extension, $G$ is a split reductive group over $K$ and we need to 
study algebraic (in)dependence in the ring $K_0(\uRep_K G).$
We assume in addition that $G$ is connected and fix a $K$-split maximal torus
$S\subset G$. Then restriction induces an injective ring homomorphism

\[
K_0(\uRep_K G)\hookrightarrow K_0(\uRep_K S)
\]

\cite{SerreRepofred}, 3.4 and the structure of $K_0(\uRep_K S)$ is particularly
simple, c.f. (\ref{fünf}) for the more general case of a diagonalizable group. It then becomes important to determine the character group of $S$
which in the present case (the topological group $T$ being arbitrary) we do not know how to do.
The most satisfactory solution of our initial problem can be obtained in the case that $T=<F>$ is topologically cyclic. Though very special, this case will be of interest 
for us because it covers the case of the Galois group of a finite field.
So let 

\[
\phi: T=<F>\longrightarrow\Gl(V)(K)\simeq\Gl_m(K)
\]

be as in (\ref{sternchen}) (so $m=\sum_{i=1}^n\ddim_K V_i$). By our preliminary reductions (\ref{eins}) and
(\ref{vier}) and conjugating $\phi(F)$ suitably we can assume that

\begin{equation}\label{tilde}
\phi(F)=\diag(\alpha_1,\ldots ,\alpha_m)\in\D_m(K)
\end{equation}

for some $\alpha_i\in K^*$, where $\D_m\subset\Gl_m$ is the torus
of diagonal matrices. Then $G$ as in (\ref{drei}) is the
smallest algebraic subgroup $G\subset\Gl_m$ defined over $K$
such that $\phi(F)\in G(K)$ and $D:=G\subset\D_m$
is a closed subgroup, hence diagonalizable by \cite{Borel}, proposition 8.4.
The character group $X(D)$ of $D$ can be identified with the subgroup of $K^*$ generated by the eigenvalues of $\phi(F)$ as
follows from \cite{Borel}, 8.2, corollary:

\begin{equation}\label{aleph}
X(D)\simeq<\alpha_1,\ldots, \alpha_m>\subset K^*.
\end{equation}

We recall the structure of $K_0(\uRep_K D)$ from \cite{SerreRepofred}, 3.4.
Giving $V\in\uRep_K D$ is equivalent to giving a ``grading of type $X(D)$'' 
of the finite dimensional $K$-vectorspace $V$, i.e. $V=\oplus_{\chi\in X(D)} V_{\chi}$,
where $V_{\chi}\subset V$ is the $\chi$-eigenspace and 

\begin{equation}\label{fünf}
\ch: K_0(\uRep_K D)\stackrel{\sim}{\longrightarrow}\Z[X(D)]\; , \; [V]\mapsto\sum_{\chi} \dim_K(V_{\chi})e^{\chi}
\end{equation}

is a ring isomorphism. Here $\Z[X(D)]$ is the group algebra over $X(D)$, i.e.
$\Z[X(D)]$ is $\Z$-free on the set $\{e^{\chi}:\chi\in X(D)\}$ and
$e^{\chi}e^{\chi'}=e^{\chi+\chi'}$.
Decomposing the finitely generated abelian group $X(D)$ as

\[
X(D)=X(D)^{tor}\oplus_{i=1}^d \Z e_i,
\]

where $X(D)^{tor}\subset X(D)$ is the torsion subgroup and the $e_i\in X(D)$
are suitably chosen, the natural map

\[
\Z[X(D)^{tor}]\otimes_{\Z}\Z[\oplus_{i=1}^d\Z e_i]\stackrel{\sim}{\longrightarrow}
\Z[X(D)]\; , \; e^t\otimes e^f\mapsto e^{t+f}
\]

is an isomorphism and

\[
\Z[\oplus_{i=1}^d\Z e_i]\stackrel{\sim}{\longrightarrow}\Z[T_1,\ldots ,T_d,T_1^{-1},\ldots, T_d^{-1}]\; , \; e^{\sum a_ie_i}\mapsto T_1^{a_1}\ldots T_d^{a_d}
\]

is an isomorphism to the ring of Laurant polynomials over $\Z$ in $d=$rk$(X(D))$=dim($D$) variables. Since
$X(D)^{tor}$ is finite, $A:=\Z[X(D)^{tor}]$
is a finite flat $\Z$-algebra and summing up we have

\begin{equation}\label{kreuz}
K_0(\uRep_K D)\simeq A[T_1,\ldots ,T_d,T_1^{-1},\ldots, T_d^{-1}].
\end{equation}

Note that, if $X(D)$ is as in (\ref{aleph}), then $X(D)^{tor}\subset K^*$
is a finite subgroup, hence cyclic of order $N$, say. Then we have

\begin{equation}\label{tor}
A=\Z[\Z/N]\simeq\Z[S]/(S^N-1).
\end{equation}

So our initial problem is finally reduced to deciding whether the

\[
f_i:=\ch([V_i])\in A[T_1,\ldots ,T_d,T_1^{-1},\ldots, T_d^{-1}],
\]

i.e. the characters of the given representations $V_i$, are algebraically independent.
For this problem one can prove a Jacobian criterion as follows:\\
Let $A$ be a finite flat $\Z$-algebra. We will prove
a criterion for given 

\[
f_1,\ldots, f_n\in A[T_1,\ldots, T_m,T_1^{-1}, \ldots, T_m^{-1}]
\]
 
to be algebraically independent in terms of the Jacobian

\[
J:=\left( \frac{\partial f_i}{\partial T_j}\right)_{1\le i \le  n,\; 1\le j \le m}
\]

which is an $n$ by $m$ matrix with entries from $A[T_1,\ldots, T_m,T_1^{-1}, \ldots, T_m^{-1}]$.\\
For a homomorphism $\phi: A\rightarrow\oQ$ we denote by

\begin{eqnarray*}
\phi\otimes 1 : A[\uT,\uT^{-1}]:=A[T_1,\ldots, T_m,T_1^{-1}, \ldots, T_m^{-1}]=
A\otimes\Z[\uT,\uT^{-1}]\rightarrow \\
\rightarrow\oQ[\uT,\uT^{-1}]:=\oQ[T_1,\ldots, T_m,T_1^{-1}, \ldots, T_m^{-1}]=
\oQ\otimes\Z[\uT,\uT^{-1}]
\end{eqnarray*}

the base change of $\phi$. Then $(\phi\otimes 1)(J)$ is a matrix with coefficients
in 

\[
\oQ[\uT,\uT^{-1}]\subset Q:=\oQ(T_1,\ldots, T_m),
\]

where $Q$ is the field of rational functions in $m$ variables over $\oQ$.
For a homomorphism $\psi: \oQ[\uT,\uT^{-1}]\longrightarrow K$
with $K$ some field $(\psi\circ(\phi\otimes 1))(J)$
is a matrix with coefficients in $K$.
Finally, by $\rk_L(M)$ we will denote the rank of a matrix $M$
having coefficients in the field $L$.
The Jacobi criterion for algebraic independence
mentioned above is the following.

\begin{lemma}\label{jacobigross}

Keeping the above notations we have:\\

\begin{tabular}{cl}
1) & The following are equivalent:\\
  & i) The $f_i\in A[\uT,\uT^{-1}]$ are algebraically dependent.\\
  & ii) For all $\phi: A\rightarrow\oQ$ we have
\end{tabular}

\[
\rk_Q((\phi\otimes 1)(J))<n.
\]

\begin{tabular}{cl}

2) & The following are equivalent:\\
   & i) The $f_i\in A[\uT,\uT^{-1}]$ are algebraically independent.\\
   & ii) There are homomorphisms $\phi: A\rightarrow\oQ$ and $\psi: \oQ[\uT,\uT^{-1}]\rightarrow K$\\
   & ($K$ some field) such that
\end{tabular}

\[
\rk_K(\psi\circ(\phi\otimes 1)(J))=n.
\]

\end{lemma}

The proof of lemma \ref{jacobigross} consists in a reduction, which we will omit, 
to the case when $A$ is a field of characteristic zero. In this case the lemma is \cite{Borel}, chapter AG, theorem 17.3.
We conclude this subsection by giving two applications.

\begin{theorem}\label{abelian surfaces}

Let $k$ be a finite field and let $E_1,E_2/k$ be non-isogeneous ordinary elliptic curves. Let $E'_2/k$ be the quadratic twist of $E_2$, and consider the
abelian surfaces $A_1:=E_1\times E_2$ and $A_2:=E_1\times E'_2$ over $k$
and the associated Galois representations
$V_i:=H_c^1(\overline{A_i})\in\Rep_{G_k}\Q_l$.
Let $k\subset L$ be the unique quadratic extension of $k$ inside $\ok$. Then:\\
i) $\Res^{G_k}_{G_L}(V_1)\simeq\Res^{G_k}_{G_L}(V_2)$.\\
ii) The classes $[V_1]$ and $[V_2]$ are algebraically independent in $K_0(\Rep_{G_k}\Q_l)$.\\
In particular, the classes $[A_1]$ and $[A_2]$ are algebraically independent in
$K_0(\Var_k)$.

\end{theorem}

\begin{rem} The representations $V_1$ and $V_2$ are algebraically independent but
not geometrically algebraically independent. In fact, after restriction to
the subgroup of index two of $G_k$ their classes become algebraically
dependent and in fact equal. We also see that $x:=[A_1]-[A_2]\in K_0(\Var_k)$
generates a polynomial ring over $\Z$ inside $K_0(\Var_k)$ even though
$x$ lies in the kernel of the base change homomorphism $K_0(\Var_k)\longrightarrow K_0(\Var_L)$. This shows that $K_0(\Var_k)$ encodes fine arithmetic
invariants of varieties over $k$.

\end{rem}

\begin{proof} Part i) is clear because $A_1\times_k L\simeq A_2\times_k L$ and the last assertion follows from part ii) and corollary
\ref{H1suffices}. We prove ii):
Fix $\pi_i\in\oQ^*$ a Weil number attached to $E_i/k$ and denote by
$\opi_i$ the conjugate of $\pi_i$. Note that $\Q(\pi_1)$ and $\Q(\pi_2)$  
are distinct imaginary quadratic fields, hence

\begin{equation}\label{mann}
\Q(\pi_1)\cap\Q(\pi_2)=\Q.
\end{equation}

The eigenvalues of the geometric Frobenius $F_k$ acting on the $H_c^1$ of
$E_1$,$E_2$ and $E'_2$ are (respectively) $\{\pi_1,\opi_1\}$, $\{\pi_2,\opi_2\}$
and $\{-\pi_2,-\opi_2\}$. Fixing a finite extension $\Q_l\subset K$ containing $\pi_1$
and $\pi_2$ we get in the notation of (\ref{tilde}):
\[
\phi(F_k)=\diag(\pi_1,\opi_1,\pi_2,\opi_2,\pi_1,\opi_1,-\pi_2,-\opi_2)\in\D_8(K)
\]
and
\[
X(D)\simeq<\pi_1,\opi_1,\pi_2,\opi_2,-1>\subset K^*.
\]

Using (\ref{mann}) and the fact that $E_1$
and $E_2$ are ordinary one checks that 
the structure of $X(D)$ is as follows.

\begin{equation}\label{structure}
X(D)=\{\pm 1\}\oplus\pi_1^{\Z}\oplus(\opi_1)^{\Z}\oplus\pi_2^{\Z}.
\end{equation}

From (\ref{structure}) we get as in (\ref{kreuz}) together with (\ref{tor})

\begin{equation}\label{frage}
K_0(\uRep_K D)\simeq (\Z[S]/(S^2-1))[T_1,T_1^{-1},\ldots,T_3,T_3^{-1}]
\end{equation}

and compute the characters $f_i:=\ch([V_i])$ as $f_1=e^{\pi_1}+e^{\opi_1}+e^{\pi_2}+e^{\opi_2}$
and $f_2=e^{\pi_1}+e^{\opi_1}+e^{-\pi_2}+e^{-\opi_2}$.
We have $\opi_2=\pi_1\opi_1\pi_2^{-1}$
and using this, (\ref{structure}) and (\ref{frage}) we get (denoting abusively the
images of $[V_i]$ under (\ref{frage}) by $f_i$ again)

\[
f_1=T_1+T_2+T_3+T_1T_2T_3^{-1}\mbox{ and}
\]
\[
f_2=T_1+T_2+ST_3+ST_1T_2T_3^{-1}.
\]

To show that $f_1$ and $f_2$ are algebraically independent we compute the Jacobian of $f_1$ and $f_2$

\[
J=\left(\begin{array}{ccc}
         1+T_2T_3^{-1} & 1+T_1T_3^{-1} & 1-T_1T_2T_3^{-2} \\
         1+ST_2T_3^{-1} & 1+ST_1T_3^{-1} & S-ST_1T_2T_3^{-2}
\end{array}\right)
\]

which we specialize by putting $S=-1$, $T_1=-1$, $T_2=T_3=1$ to obtain

\[
\left(\begin{array}{ccc}
         2 & 0 & 2 \\
         0 & 2 & -2
\end{array}\right)
\]

which has rank two, hence $f_1$ and $f_2$ are algebraically independent
by lemma \ref{jacobigross}, as was to be shown.

\end{proof}

Our next application of the methods developed so far is the following.

\begin{theorem}\label{inftycurvesfinitefield}
Let $k$ be a finite field. Then there is a sequence of proper, smooth and geometrically connected
curves $X_i/k$ ($i\ge 1$) such that the classes $[X_i]\in K_0(\Var_k)$ 
are algebraically independent.
\end{theorem}

\begin{theorem}\label{ellovernumber}
Let $k$ be a number field and $\{E_i\}_{i\in I}$ a set of elliptic curves 
over $k$ such that the $E_i$ are pairwise non-isogeneous and satisfy
$\End_{\overline{k}}(E_i)=\Z$. Then the classes $[E_i]\in K_0(\Var_k)$ are
algebraically independent.
\end{theorem}

Note that over any number field there are infinite sets of elliptic curves satisfying the assumptions of theorem \ref{ellovernumber}.\\
To prove these results we will need the following sufficient criterion for 
geometric algebraic independence in terms of Frobenius eigenvalues.

\begin{lemma}\label{geoindepofprocyclic}
Let $\hat{\Z}$ denote the pro-finite completion of $\Z$ and let $K/\Q_l$
be a finite extension and let $V_1,\ldots ,V_n\in\Rep_{\hat{\Z}} K$ be given. For $1\leq i\leq n$ put $m_i:=\ddim_K V_i$
and let
\[
\{ \lambda_{i,j}\}_{1\le j\le m_i}\subset\overline K^*
\]
be the eigenvalues of $F:=1$ acting on $V_i$.
For any $1\le k\le n$ consider the finite dimensional $\Q$-vectorspace
\begin{equation}\label{A_k}
A_k:=<\lambda_{i,j}>_{1\leq i\leq k,\; 1\le j\le m_i}\otimes_{\Z}\Q
\end{equation}

and set $A_0:=0$. If

\begin{equation}\label{dimA_k}
\ddim_{\Q} A_k>\ddim_{\Q} A_{k-1}\mbox{ for k}=1,\ldots, n
\end{equation}

then the $V_1,\ldots, V_n\in \Rep_{\hat{\Z}} K$ are geometrically 
algebraically independent.

\end{lemma}

For the application of this lemma in the proof of theorem \ref{ellovernumber}
we make its assumptions more explicit in the case when the $V_i$ are (dual to) 
the Tate-modules of elliptic curves.

\begin{lemma}\label{assumptionsell}
Let $\F_q$ be a finite field and identify $G_{\F_q}=\hat{\Z}$ using the
geometric Frobenius. Let $E_1,\ldots, E_n/\F_q$ be elliptic curves, assume
that $E_1$ is ordinary, put $V_i:=H^1_c(\overline{E_i})\in\Rep_{G_{\F_q}}(\Q_l)$
and define the $A_i$ as in lemma \ref{geoindepofprocyclic}. Then the assumptions
of this lemma are satisfied if and only if $\dim_{\Q} A_n=n+1$.
\end{lemma}

\begin{proof} Since $E_1$ is ordinary we have $\dim_{\Q} A_1=2$. For 
$k=2,\ldots, n$ we have $\dim_{\Q} A_k\leq\dim_{\Q} A_{k-1}+1$ because
$\lambda_{k,1}\lambda_{k,2}=q\in A_{k-1}$. The assertion is obvious now.
\end{proof}

\begin{proofof} lemma \ref{geoindepofprocyclic}. In the notation of (\ref{tilde}) we have

\[
\phi(F)=\diag(\lambda_{1,1},\ldots,\lambda_{1,m_1},\ldots,\lambda_{n,1},\ldots, 
\lambda_{n,m_n})\in\D_{\sum_{i=1}^n m_i}(K).
\]

As a first step we will prove that the characters $f_i:=\ch([V_i])\in K_0(\uRep_K D)$
are algebraically independent. Here $D\subset\D_{\sum_{i=1}^n m_i}$ is the smallest
algebraic subgroup with $\phi(F)\in D(K)$. After suitably enlarging $K$
(legitimate by (\ref{vier})) we identify (as in (\ref{aleph}))
the character group $X(D)$ of $D$

\begin{equation}\label{XD=}
X(D)\simeq<\lambda_{i,j}>_{1\leq i\leq n,\; 1\leq j\leq m_i}\subset K^*
\end{equation}

and choose a decomposition  

\begin{equation}\label{structXD}
X(D)=<\zeta>\oplus_{i=1}^N\pi_i^{\Z}
\end{equation}

with $\zeta\in K^*$ a root of unity of order $M$, say, and 
$\pi_1,\ldots, \pi_N\in K^*$ a $\Z$-basis of the free part of $X(D)$.
We order the $\pi_i$'s such that their images $\pi_i\otimes 1$ in $X(D)\otimes\Q
=A_n$ satisfy the following: For all $1\leq k \leq n\;\; \pi_1\otimes 1,\ldots, \pi_{d_k}\otimes 1$ is a basis of $A_k$.
Here $d_k:=\dim_{\Q} A_k$ and (\ref{dimA_k}) gives

\begin{equation}\label{dks}
d_0=0<d_1<\ldots<d_n=N.
\end{equation}

Note that

\begin{eqnarray}
& & N=\ddim(D)=\ddim_{\Q}(X(D)\otimes\Q))= \label{Ngen}\\
& & \stackrel{(\ref{XD=})}{=}\ddim_{\Q}(<\lambda_{i,j}>_{1\leq i\leq n,\; 1\leq j \leq m_i}\otimes\Q)= \nonumber\\
& & \stackrel{(\ref{A_k})}{=}\ddim_{\Q}(A_n)\stackrel{(\ref{dimA_k})}{\ge}n.\nonumber
\end{eqnarray}

For any $1\leq i\leq n$ and $1\leq j\leq m_i$ we decompose, according to 
(\ref{structXD}), 

\begin{equation}\label{base}
\lambda_{i,j}=\zeta^{a_{i,j}}\pi_1^{a_{1,i,j}}\ldots\pi_N^{a_{N,i,j}}
\end{equation}

for suitable $a_{i,j},a_{k,i,j}\in\Z$. Using (\ref{kreuz})

\[
K_0(\uRep_K D)\simeq A[T_1,T_1^{-1},\ldots,T_N,T_N^{-1}]
\]

where, according to (\ref{tor}), we have $A\simeq\Z[S]/(S^M-1)$.
The characters under consideration are given as follows (taking
(\ref{frage}) as an identification)

\begin{eqnarray}
& & f_i=\ch([V_i])=\sum_{j=1}^{m_i}\lambda_{i,j}\stackrel{(\ref{base})}{=}
\sum_{j=1}^{m_i}\zeta^{a_{i,j}}\pi_1^{a_{1,i,j}}\ldots\pi_N^{a_{N,i,j}}= \label{expansion}\\
& & =\sum_{j=1}^{m_i}S^{a_{i,j}} T_1^{a_{1,i,j}}\ldots T _N^{a_{N,i,j}},\; 1\leq i\leq n.\nonumber
\end{eqnarray}

We define $\phi:A\longrightarrow\oQ$ by $S\mapsto 1$ and

\[
\psi:\oQ[T_1,T_1^{-1},\ldots,T_N,T_N^{-1}]\hookrightarrow Q:=\oQ(T_1,\ldots, T_N)
\]

to be the natural inclusion.
Now we show that the $[V_i]\in K_0(\Rep_{\hat{\Z}} K)$ are algebraically independent.
Using lemma \ref{jacobigross} it suffices to show
that $\psi\circ(\phi\otimes 1)(J)$ has rank $n$ over $Q$. Here $J$ is
the Jacobi matrix $J=\left(\frac{\partial f_i}{\partial T_j}\right)_{1\leq i\leq n,\; 1\leq j\leq N}$.
Since $n\leq N$ by (\ref{Ngen}) this is equivalent to the vectors

\[
v_i:=\left(\psi\circ\left(\phi\otimes 1\right)\left(\frac{\partial f_i}{\partial T_j}\right)\right)_{1\leq j\leq N}\; , \; 1\leq i\leq n
\]

being linearly independent over $Q$. Let us abbreviate $x_{i,j}:=\psi\circ\left(\phi\otimes 1\right)\left(\frac{\partial f_i}{\partial T_j}\right)\in Q$, 
hence $v_i=(x_{i,j})_{1\leq j\leq N}$.
We first prove the following claim concerning the shape of the matrix
$(x_{i,j})$.\\
{\em Claim:} For every $1\leq i\leq n$ there is some $1\leq j_i\leq N$
such that

\begin{eqnarray}
x_{i',j_i} & = & 0\mbox{ for }i'<i\mbox{ and}\label{shape}\\
x_{i,j_i} & \neq & 0.\nonumber
\end{eqnarray}

Indeed, given $1\leq i\leq n$ let $j_i$ be any integer satisfying 
$d_{i-1}<j_i\leq d_i$ which exists by (\ref{dks}). This choice 
guarantees that, in the notation of (\ref{expansion}), we have

\[
a_{j_i,i,l}\geq 1\mbox{ for some }1\leq l\leq m_i\mbox{ and}
\]
\[
a_{j_i,i',l}=0 \mbox{ for }i'<i\mbox{ and all }l,
\]
the claim follows.
Put a little less formally: In (\ref{expansion}) there is, for 
every $i$, some variable $T_{j_i}$ which occurs in $f_i$
but does not occur in any $f_{i'}$ with $i'<i$. This is where (\ref{dimA_k}) is used.
Now the linear independence of $v_1,\ldots, v_n$ over $Q$ follows.
Assume by contradiction that we are given $\alpha_1,\ldots,\alpha_k\in Q$, 
$\alpha_k\neq 0$, such that $\sum_{i=1}^k\alpha_i v_i=0$.
Taking the $j_k$-component of this relation we get
$\sum_{i=1}^k\alpha_i x_{i,j_k}=0\mbox{ in }Q$
which by (\ref{shape}) simplifies to
$\alpha_k x_{k,j_k}=0\mbox{ with } x_{k,j_k}\neq 0$,
hence $\alpha_k=0$, a contradiction.
At this point the proof of the fact that the classes $[V_i]\in K_0(\Rep_{\hat{\Z}} K)$ are algebraically independent is complete.\\
To show that they are in fact geometrically algebraically independent, let an open subgroup
of finite index $T'\subset\hat{\Z}$ be give. There is an integer $L\ge 1$
such that $T'=L\hat{\Z}$ and the set of eigenvalues of the topological generator
$L\in T'$ acting on $\Res^{\hat{\Z}}_{T'}(V_i)$ is $\{ \lambda_{i,j}'\}_{i,j} = \{ \lambda_{i,j}^L \}_{i,j}$. So the homomorphism

\[
<\lambda_{i,j}>_{1\leq i\leq k,1 \leq j \leq m_i}\longrightarrow
<\lambda_{i,j}'>_{1\leq i\leq k,1 \leq j \leq m_i}\;,\; x\mapsto x^L
\]

is surjective with finite kernel and induces an isomorphism

\[
A_k\stackrel{\sim}{\longrightarrow} <\lambda_{i,j}'>_{1\leq i\leq k,1 \leq j \leq m_i}\otimes_{\Z}\Q
\]

for all $1\leq k \leq n$. Thus the tuple $(T', \Res^{\hat{\Z}}_{T'}(V_i))$
satisfies the analogue of condition (\ref{dimA_k}) and by
what has already been proved we see that the $[\Res^{\hat{\Z}}_{T'}(V_i)]\in K_0(\Rep_{T'}K)$ are
algebraically independent. As $T'$ was arbitrary this shows that the $V_i$
are in fact geometrically algebraically independent and concludes the proof.
\end{proofof}

Before engaging into the slightly involved proof of theorem \ref{inftycurvesfinitefield}
we remark that if in this result one replaces curves by abelian varieties 
one can give a quick proof as follows: The subgroup $\Gamma\subset\overline{\Q}^*$
generated by $q$-Weil numbers of weight one has infinite rank. Choosing a suitable
sequence of such Weil numbers and corresponding abelian varieties (which
exist by T. Honda) one obtains the result using lemma \ref{geoindepofprocyclic}
and corollary \ref{H1suffices}. Even though $\rk(\Gamma)=\infty$ is well-known
we could not find it stated explicitly in the literature so we record 
this result (which will be clear from our proof of theorem \ref{inftycurvesfinitefield}) for reference.

\begin{cor} Let $q$ be a prime power. Then the subgroup of
$\overline{\Q}^*$ generated by $q$-Weil numbers of weight one has infinite
rank.
\end{cor}

We now start the proof of  theorem \ref{inftycurvesfinitefield}.
In fact, we will give a family of curves satisfying the conclusion of
theorem \ref{inftycurvesfinitefield} explicitly.\\
Let $\F_p\subset k$ be the prime field of $k$ and for any $a\ge 1$
let $Y_a/\F_p$ be the proper smooth curve which is birational to
the affine plane curve given by

\begin{eqnarray*}
y^p-y & = & x^{p^a-1}.
\end{eqnarray*}

The $Y_a$ are geometrically connected of genus $\frac{1}{2}(p^a-2)(p-1)$.
The arithmetic of these curves has been studied by H. Davenport and 
H. Hasse \cite{DavenportHasse}. Unfortunately, we cannot assert that 
the sequence $\{ Y_a \}_{a\ge 1}$ satisfies the conclusion of theorem \ref{inftycurvesfinitefield}
but a suitable subsequence will do.
Let $T=T(p)$ be the integer determined by lemma
\ref{phigrows} below and let $\{ n_i \}_{i\ge 1}$ be the sequence of prime numbers
greater than $T$ ordered increasingly. Then the

\[
X_i:=Y_{n_i}\otimes_{\F_p} k \; , \; i\ge 1
\]

are proper, smooth and geometrically connected curves over $k$ 
and the following result makes the existence statement of theorem
\ref{inftycurvesfinitefield} explicit.

\begin{theorem}\label{DHcurvesdo} 
The classes $[X_i]\in K_0(\Var_k)$ ($i\ge 1$) are algebraically independent.
\end{theorem}

For any $i\ge 1$ let us write $V_i:=H^1_c(\oY_{n_i})\in\Rep_{G_{\F_p}}\Q_l$.
Using corollary \ref{H1suffices}, theorem \ref{DHcurvesdo} results from the following.

\begin{theorem}\label{geoindepofH1}
The classes $V_i\in\Rep_{G_{\F_p}}\Q_l$ ($i\ge 1$) are
geometrically algebraically independent.
\end{theorem}

This will be proved using lemma \ref{geoindepofprocyclic} for which we need to review some results 
concerning the Weil numbers attached to the Jacobians of the curves $Y_a/\F_p$.\\
We fix a valuation $v$ of $\oQ$ lying over $p$ and normalized by $v(p)=1$.
We will denote by the same letter $v$ the restriction of $v$ to any
subfield $K\subset\oQ$.\\
Let $F_p$ be the geometric Frobenius of $\F_p$. Then the eigenvalues of
$F^a_p$  acting on $H^1_c(\oY_a)$ are given by generalized Jacobi sums

\begin{equation}\label{jacsum}
\tau_j(\chi^t)^{(a)}:=-\sum_{u\in\F_{p^a}^*}\chi^t(u)\exp\left( \frac{2\pi i}{p}
j\tr (u) \right)\; , \; 1\le j\le p-1,\; 1\le t\le p^a-2,
\end{equation}

where $\chi:\F_{p^a}^*\hookrightarrow\C^*$ is a certain faithful character
of $\F_{p^a}^*$ (a Teichm\"uller lift depending on the choice of $v$, see \cite{Manin}, \S4) and $\tr: \F_{p^a}\longrightarrow \F_p$ denotes the trace,
\cite{Yamada}, (2).\\
We will abbreviate $\tau(\chi^t)^{(a)}:=\tau_1(\chi^t)^{(a)}$ in the 
following.
Let $1\le t \le p^a-2$ be given and let 

\[
t=\sum_{i=0}^{a-1}j_ip^i\; , \; 0\le j_i\le p-1
\]

be the $p$-adic expansion of $t$ and put $\sigma(t):=\sum_{i=0}^{a-1}j_i$.
Then \cite{Manin}, \S 4 gives the valuations of the generalized Jacobi sums as

\begin{equation}\label{randale}
v(\tau(\chi^t)^{(a)})=\frac{\sigma(t)}{p-1}.
\end{equation}

We will need the following lower bound on Euler's phi-function, 
$\phi(n):=|(\Z/n)^*|$ for $n\ge 1$.

\begin{lemma}\label{phigrows}

Let $p$ be a prime number. Then there is an integer $T=T(p)$ such that
\[
\frac{\phi(p^n-1)}{n}>\phi(p-1)\mbox{ for all } n\ge T.
\]
\end{lemma}

\begin{proof} This follows from the inequality, valid for any $m\ge 1$,

\[
m-1\ge\phi(m)\ge\frac{m}{C\log\log m}
\]

with a suitable positive constant $C$. The first inequality is 
trivial and the second is \cite{rosserschoenfeld}, theorem 15.
\end{proof}

\begin{proofof} theorem \ref{geoindepofH1}.
It suffices to show that for any given $M\ge 1$ the $V_1,\ldots, V_M
\in\Rep_{G_{\F_p}}\Q_l$ are geometrically algebraically independent. So we fix some $M\ge 1$. For any
$1\le i\le M$ let $\{ \lambda_{i,j} \}_{1\leq j\leq (p^{n_i}-2)(p-1)}\subset\oQ^*$
be the eigenvalues of $F_p$ acting on $V_i$. Note that

\begin{equation}\label{f5}
\{ \lambda_{i,j}^{n_i} \}_{1\leq j \leq (p^{n_i}-2)(p-1)}=
\{ \tau_j(\chi^t)^{(n_i)} \}_{1\leq j \leq p-1,\; 1\leq t \leq p^{n_i}-2}
\end{equation}

by (\ref{jacsum}).
We define $B_0:=\{ 1 \} \subset\oQ^*$ and for $1\leq k\leq M$ we define
$B_k:=<\lambda_{i,j}>_{1\leq i\leq k,\mbox{{\tiny all}} j} \subset\oQ^*$
to be the subgroup of $\oQ^*$ generated by the $F_p$-eigenvalues occurring among
the $V_1,\ldots, V_k$ and 
$A_k:=B_k\otimes_{\Z}\Q\; , \; d_k:=\ddim_{\Q} A_k$.
Using lemma \ref{geoindepofprocyclic} the proof will be complete if we establish the following.

\begin{equation}\label{weilnogrow}
\mbox{ We have }d_k>d_{k-1}\mbox{ for }k=1,\ldots, M.
\end{equation}

As $B_{k-1}\subset B_k$ we always have $d_k\geq d_{k-1}$ and we need
to show that this inequality is strict.
For any $n\ge 1$ we denote by $K_n\subset\oQ$ the field of $n$-th
roots of unity and we fix isomorphisms
$\Gal(K_n/\Q)\simeq(\Z/n)^*$
using a compatible system of roots of unity in $\oQ^*$. Recall that
$K_n\cap K_m=K_{\ggcd(n,m)}$ and $K_n K_m=K_{\lcm(n,m)}$ where
$\ggcd$ and $\lcm$ mean greatest common divisor and least common multiple, 
respectively.
From (\ref{f5}) and (\ref{jacsum}) we see that

\begin{equation}\label{moses}
\lambda_{i,j}^{n_i}\in K_{p^{n_i}-1}K_p=K_{p(p^{n_i}-1)},
\end{equation}

hence, introducing $N_k:=\prod_{i=1}^k n_i$ for $0\leq k \leq M$, we have that

\begin{equation}\label{krabbe}
B_k^{N_k}:=\{ x^{N_k}:x\in B_k \} \subset  K_{p(p^{N_k}-1)}\mbox{ for } k=0,\ldots ,M.
\end{equation}

We now prove (\ref{weilnogrow}) by contradiction, assuming that we are given some 
$1\le k \le M$ such that

\begin{equation}\label{false}
d_k=d_{k-1}.
\end{equation}

As $k$ will be fixed from now on, we put $n:=n_k$ and $N:=N_{k-1}$
to ease the reading.
The cyclotomic fields to play a role in the following are organized in the 
following diagram.

\begin{equation}\label{106}
\xymatrix{ & K_{p(p^n-1)} \ar@{-}[ld] \ar@{-}[rrrd] & & K_{p(p^N-1)} \ar@{-}[ld]
\ar@{-}[rd] & & \\
K_{p^n-1} \ar@{-}[rd] & & K_{p^N-1} \ar@{-}[ld] & & K_p \ar@{-}[llldd]\\
 & K_{p-1} \ar@{-}[d] & & & \\
 & \Q & & & \\ }
\end{equation}

Observe that by construction $\ggcd(n,N)=1$ ( $N$ is a product of 
primes different from the prime $n$) and hence $\ggcd(p^n-1,p^N-1)=p-1$
which implies as recalled above that $K_{p^n-1}\cap K_{p^N-1}=K_{p-1}$
as asserted by (\ref{106}). The remaining assertions implicit in (\ref{106})
are checked similarly.
According to (\ref{jacsum}) there are $\lambda_1$ and $\lambda_2$ 
which are eigenvalues of $F_p$ acting on $V_k$ (hence $\lambda_i\in B_k$ 
by the definition of $B_k$) and such that

\begin{eqnarray}
\lambda_1^n & = & \tau(\chi^1)^{(n)} \label{piggy}\\
\lambda_2^n & = & \tau(\chi^{p^n-2})^{(n)}, \nonumber 
\end{eqnarray}

hence $\lambda_i^n\in K_{p(p^n-1)}$ by (\ref{jacsum}) and by \cite{Yamada}, (14)
we have 

\begin{equation}\label{fulton}
\tau_i:=\lambda_i^{n(p-1)}\in K_{p^n-1}\; , \; i=1,2.
\end{equation}

The idea now is the following.
We will exhibit two places $v_1$ and $v_2$ of $K_{p^n-1}$ lying above the same place
of $K_{p-1}$ such that a suitable linear combination $x$ of $\tau_1$ and 
$\tau_2$ has different valuations at $v_1$ and $v_2$ and using (\ref{false}) 
we will see that we may chose $x$ to lie in $K_{p-1}$ which will give the
desired contradiction.
Here are the details:\\
The inclusion $K_{p-1}\subset K_{p^n-1}$ corresponds to a surjection of 
Galois groups

\[ \xymatrix{ \pi:(\Z/p^n-1)^* \ar@{>>}[r] & (\Z/p-1)^*. } \]

The decomposition group of the place $v_1:=v$ of $K_{p^n-1}$ is generated
by the residue of $p$ (since $p$ is unramified in $K_{p^n-1}$) and $\pi$
factors over a surjection

\begin{equation}\label{52}
\xymatrix{ \opi:(\Z/p^n-1)^*/<\overline{p}> \ar@{>>}[r] & (\Z/p-1)^*.}
\end{equation}

The left hand side of (\ref{52}) acts simply transitively by conjugation
on the set of places of $K_{p^n-1}$ lying above $p$ and this will be
used to construct a second such place $v_2$.
The order of the left hand side of (\ref{52}) is $\phi(p^n-1)/n> \phi(p-1)=|(\Z/p-1)^*|$ by lemma \ref{phigrows} and our choice that $n\ge T(p)$
so we may choose 

\begin{equation}\label{tanz}
1\neq\gamma\in\kker(\opi).
\end{equation}

We denote by $v_2$ the place of $K_{p^n-1}$ conjugate by $\gamma^{-1}$ to $v_1$, i.e. $v_2(x)=v_1(\gamma . x)=v(\gamma . x)$ for all $x\in K_{p^n-1}^*$.
Now we compute some valuations.

\[
v_1(\tau_1)\stackrel{(\ref{fulton})}{=}(p-1)v(\lambda_1^n)\stackrel{(\ref{piggy}),(\ref{randale})}{=}(p-1)\frac{\sigma(1)}{p-1}=1\mbox{ and }
\]
\[
v_1(\tau_2)\stackrel{(\ref{fulton})}{=}(p-1)v(\lambda_2^n)\stackrel{(\ref{piggy}),(\ref{randale})}{=}(p-1)\frac{\sigma(p^n-2)}{p-1}=n(p-1)-1,
\]

using for the computation of $\sigma(p^n-2)$ that

\[
p^n-2=p^n-1-1=\sum_{i=0}^{n-1}(p-1)p^i-1=(p-2)+\sum_{i=1}^{n-1}(p-1)p^i.
\]

Let the $p$-adic expansion of $\gamma$ be given by

\begin{equation}\label{tinneff}
\gamma=\sum_{i=0}^{n-1}\alpha_i p^i\; , \; 0\leq\alpha_i \leq p-1.
\end{equation}

Observe that 

\[
\sigma(\gamma)=\sum_{i=0}^{n-1}\alpha_i\equiv\opi(\gamma)\stackrel{(\ref{tanz})}{\equiv}1\;(\mbox{mod }(p-1)).
\]

If we had $\sigma(\gamma)=1$ then \cite{Yamada}, lemma 3, ii) would 
imply that $\gamma\in p^{\Z}$, hence $\gamma=1$ in $(\Z/p^n-1)^*/<\overline{p}>$
contrary to our choice (\ref{tanz}). So we have

\begin{equation}\label{banane}
\sigma(\gamma)=1+a(p-1)\mbox{ for some }a\ge 1.
\end{equation}

Now we compute the valuations at $v_2$.
\begin{eqnarray*}
v_2(\tau_1)\stackrel{(\ref{fulton}),(\ref{piggy})}{=}(p-1)v_2(\tau(\chi^1)^{(n)})=(p-1)v(\gamma . \tau(\chi^1)^{(n)})= & & \\
=(p-1)v(\tau(\chi^{\gamma})^{(n)})\stackrel{(\ref{randale})}{=}(p-1)\frac{\sigma(\gamma)}{p-1}\stackrel{(\ref{banane})}{=}1+a(p-1).
\end{eqnarray*}

In the third term we consider $\gamma$ as an element of $\Gal(K_{p^n-1}/\Q)$
and the second equality follows from our definition of $v_2$. The third
equality uses $\gamma . \tau(\chi^t)^{(n)}=\tau(\chi^{\gamma t})^{(n)}$
which follows from (\ref{jacsum}). Similarly we have

\begin{eqnarray*}
v_2(\tau_2)\stackrel{(\ref{fulton}),(\ref{piggy})}{=}(p-1)v_2(\tau(\chi^{p^n-2})^{(n)})=(p-1)v(\tau(\chi^{\gamma (p^n-2)})^{(n)})= & & \\
=(p-1)v(\tau(\chi^{-\gamma})^{(n)})\stackrel{(\ref{randale})}{=}(p-1)\frac{\sigma(-\gamma)}{p-1}=\sigma(-\gamma)\stackrel{(*)}{=}n(p-1)-\sigma(\gamma) & & \\
\stackrel{(\ref{banane})}{=}n(p-1)-(1+a(p-1))=(n-a)(p-1)-1. & & 
\end{eqnarray*}

In this sequence of equalities, we would like to justify (*) which is an elementary property of
the function $\sigma$:
From (\ref{tinneff}) we get the $p$-adic expansion of $-\gamma$:

\[
-\gamma=-\sum_{i=0}^{n-1}\alpha_i p^i\stackrel{(\mbox{ mod }p^n-1)}{\equiv}
p^n-1-\sum_{i=0}^{n-1}\alpha_i p^i=\sum_{i=0}^{n-1}(p-1-\alpha_i)p^i,
\]

hence 

\[
\sigma(\gamma)+\sigma(-\gamma)=\sum_{i=0}^{n-1}\alpha_i+\sum_{i=0}^{n-1}(p-1-\alpha_i)=n(p-1),
\]

as used in (*).
To sum up, we have computed the matrix

\[
X:=\left(\begin{array}{cc}
          v_1(\tau_1) & v_1(\tau_2) \\
          v_2(\tau_1) & v_2(\tau_2) 
\end{array} \right)=
\left(\begin{array}{cc}
       1 & n(p-1)-1 \\
       1+a(p-1) & (n-a)(p-1)-1 
\end{array} \right)
\]

and compute its determinant

\begin{eqnarray}
\ddet(X)=(n-a)(p-1)-1-(1+a(p-1))(n(p-1)-1)= & & \label{nonsing}\\
= (n-a)(p-1)-1-(n(p-1)-1+an(p-1)^2-a(p-1)) =& & \nonumber\\
= -an(p-1)^2\neq 0, \nonumber
\end{eqnarray}

because $a\neq 0$ by (\ref{banane}).
If $\mu,\nu\in\Z$ are given and we put $x:=\tau_1^{\mu}\tau_2^{\nu}\in K_{p^n-1}$
we have

\[
\left(\begin{array}{c} v_1(x) \\
                       v_2(x)
      \end{array} \right)=
X\left(\begin{array}{c} \mu \\
                       \nu
      \end{array} \right).
\]

By (\ref{nonsing}) there is some integer $R\ge 1$ and $\mu,\nu\in\Z$
such that $X\left(\begin{array}{c} \mu \\ \nu \end{array} \right)=\left(\begin{array}{c} R \\ 0 \end{array} \right)$, so if we put $x:=\tau_1^{\mu}\tau_2^{\nu}$ for these particular values of $\mu$ and $\nu$ we get $x\in K_{p^n-1}$ 
with 

\begin{equation}\label{zeugel}
v_1(x)=R\ge 1\mbox{ and }v_2(x)=0.
\end{equation}

Note that for the restrictions of $v_1$ and $v_2$ to $K_{p-1}$ we have

\begin{equation}\label{spinat}
v_1|_{K_{p-1}}=v_2|_{K_{p-1}}
\end{equation}

because the restriction of $\gamma\in\Gal(K_{p^n-1}/\Q)$ to $K_{p-1}$
is trivial, i.e. $\opi(\gamma)=1$, c.f. (\ref{tanz}).
Now we can finally derive a contradiction.
By (\ref{false}) we know that $B_k/B_{k-1}$ is torsion. As $\lambda_1,\lambda_2
\in B_k$ we find an integer $S\ge 1$ with 

\begin{equation}\label{DING}
\lambda_i^S\in B_{k-1}^N\stackrel{(\ref{krabbe})}{\subset} K_{p(p^N-1)}\; , \; i=1,2
\end{equation}

(recall that $N=N_{k-1}$) and furthermore 
$x^S=(\tau_1^{\mu}\tau_2^{\nu})^S\stackrel{(\ref{fulton})}{=} (\lambda_1^{\mu}\lambda_2^{\nu})^{n(p-1)S}$,
so $x^S\in K_{p^n-1}$ by (\ref{fulton}) and $x^S\in K_{p(p^N-1)}$ by
(\ref{DING}) from which we conclude that

\begin{equation}\label{rahm}
x^S\in K_{p^n-1}\cap K_{p(p^N-1)}\stackrel{(\ref{106})}{=} K_{p-1}.
\end{equation}

This implies
$SR\stackrel{(\ref{zeugel})}{=}v_1(x^S)\stackrel{(\ref{rahm}),(\ref{spinat})}{=}
v_2(x^S)\stackrel{(\ref{zeugel})}{=}0$,
a contradiction which concludes the proof.

\end{proofof}

\begin{proofof} theorem \ref{ellovernumber}. We can assume that the given set
of elliptic curves $\{E_i\}_{i\in I}=\{E_1,\ldots,E_n\}$ is finite. We will
show that the classes of $V_i:=H^1_c(\overline{E_i})$ in $K_0(\Rep_{G_k}\Q_l)$
are algebraically independent which is sufficient by corollary \ref{H1suffices}.
We consider the subring $R:=\Z[[V_i]]\subset K_0(\Rep_{G_k}\Q_l)$ generated 
by the $[V_i]$. For any finite place $v$ of $k$ with residue field 
$k(v)$ at which all the $E_i$ have good reduction we have a canonical
homomorphism $\phi_v:R\longrightarrow K_0(\Rep_{G_k(v)}\Q_l)$ because all
$V_i$ are unramified at $v$. We are going to construct such a place $v$
such that the $\phi_v([V_i])$ are algebraically independent. For $A:=E_1\times\cdots\times E_n$
we have $H^1_c(\overline{A})\simeq\oplus_{i=1}^n V_i$ and denote by
$G_l\subset\Gl(H^1_c(\overline{A}))(\Q_l)$ the closure of the image of Galois.
Since $\End_{\overline{k}}(E_i)=\Z$ we know that the Hodge-group $\mathrm{Hdg}(H^1(E_i))$ is isogeneous to $\mathrm{Sl}_2$. As the $E_i$ are pairwise
non-isogeneous we have that $\mathrm{Hgd}(H^1(A))$ is isogeneous to 
$\mathrm{Sl}_2^n$ and hence the Mumford-Tate group $\mathrm{MT}(H^1(A))$
has rank $n+1$. Now observe that the Mumford-Tate conjecture is known to
be true for $A/k$: For $n=1$ it is a celebrated result of J-P. Serre and 
one reduces to this case using \cite{Ribet}, lemma on page 790.
We conclude that $G_l$ has rank $n+1$. By \cite{SerreTorus} we find a finite place
$v$ of $k$ of good reduction for all $E_i$ and such that the corresponding
Frobenius torus $T_v\subset G_l^0$ is maximal, i.e. $T_v$ has rank
$n+1$. But this rank is nothing else than $\dim_{\Q} A_n$ in the notation of
lemma \ref{assumptionsell} (applied to the reductions at $v$ of the $E_i$) 
and using this lemma we conclude that the $\phi_v([V_i])$ are indeed 
algebraically independent as was to be shown.
\end{proofof}

\vspace{.5cm}

\subsection{Using a result of Skolem}\label{using a result of skolem}

In the previous subsection we used a Tannakian argument to establish
the algebraic independence of $l$-adic Galois representations given by
the cohomology of suitable varieties.
In this subsection we approach the same goal using a lemma of Skolem
on the shape of the power series expansion of rational functions together
with the rationality properties given in \ref{finite base field}. The possibility
of this was suggested by reading a note of J-P. Serre \cite{Serreskolem}.\\
Let $k$ be a finite field. We will use the motivic measure from 
proposition \ref{rationalimage} 
\[
\mu_k: K_0(\Var_k)\longrightarrow \Rc(\Q)[T]
\]

to give sufficient conditions for the classes of given varieties over
$k$ to be algebraically independent in $K_0(\Var_k)$.\\
The main observation is that, if there is a non-trivial algebraic 
relation among the motivic measures of some varieties, then, at least after
a finite extension of the base field, there will be an {\em irreducible}
relation, see \ref{irreducible}.\\
For $\nu\ge 1$ we denote by $k_{\nu}$ the unique extension of degree $\nu$
of $k$ inside some fixed algebraic closure of $k$.

\subsection{Irreducible equations} \label{irreducible}

\begin{theorem} \label{irreducibleeq}

Let $k$ be a finite field and let $X_1,\ldots,X_n/k$ $(n\ge 1)$ be separated and of finite type and assume
that $\mu_k([X_1]),\ldots,\mu_k([X_n])$ are algebraically {\em dependent} in $\Rh(\Q)[T]$. Then there are $M\ge 1$
 and an {\em irreducible} polynomial $G\in\Z[T_1,\ldots,T_n]$ such that
\[
G(\mu_{k_M}([X_1\times_k k_M]),\ldots, \mu_{k_M}([X_n\times_k k_M]))=0.
\]
\end{theorem}

\begin{rem} This result is non-trivial even for $n=1$ because 
$\Rc(\Q)[T]$ contains zero divisors. In the remark following
the proof of theorem \ref{structureofS} we give an
example showing the necessity of allowing the finite extension $k\subset k_M$
 in theorem \ref{irreducibleeq}.
\end{rem}

The starting point of the proof of theorem 
\ref{irreducibleeq} is a result due to Skolem which we now explain.
For brevity, we will call a subset $X\subset\N$ {\em good} if it
is the union of a finite set and finitely many arithmetic progressions 
for a single modulus, i.e. $X$ is good if there are $\Sigma\subset\N$ finite, 
$M\ge 1$ and $I\subset\{0,\ldots, M-1\}$ such that
\[
X=\Sigma\cup_{i\in I}(i+M\N).
\]
The collection of good subsets is stable under finite union and finite intersection.
Skolem's result is the following \cite{Serreskolem}.

\begin{prop} \label{zeroes}

Let $f\in\Q(t)$ be a rational function with power series expansion
$f=\sum_{\nu\in\Z}a_{\nu}t^{\nu}$. Then the set $\{\nu\ge 0 : a_{\nu}=0\}$
is good.

\end{prop}

We will use the following consequence of Skolem's result.

\begin{cor} \label{zeroes1}

For $f\in\Rc(\Q)$ the set $\{\nu\ge1:\partial_{\nu}(f)=0\}$ is good.

\end{cor}

Indeed, for $f\in\Rc(\Q)$ we have that dlog($f$) is a rational function
to which we apply proposition \ref{zeroes}. As the collection of good subsets is stable under
finite intersection we also have the following.

\begin{cor}\label{zeroes2}

For $F\in\Rc(\Q)[T]$ the set $\{\nu\ge1:\partial_{\nu}(F)=0 \mbox{ in } \Q[T]\}$ is good.

\end{cor}

Corollary \ref{zeroes2} follows from corollary \ref{zeroes1} applied to
the coefficients of $F$. By abuse of notation, in corollary
\ref{zeroes2}, we have written $\partial_{\nu}:\Rc(\Q)[T]\rightarrow\Q[T]$ for
the map derived from $\partial_{\nu}:\Rc(\Q)\subset\Lambda(\Q)\rightarrow\Q$ by sending $T\mapsto T$.

\begin{proofof} theorem \ref{irreducibleeq}. By assumption we have
some $0\neq H\in\Z[T_1,\ldots,T_n]$ such that
\begin{equation}\label{relation}
H(\mu_k([X_1]),\ldots,\mu_k([X_n]))=0 \mbox{ in }\Rc(\Q)[T].
\end{equation}

Decompose $H=\prod_{i=1}^N H_i$ into a product of irreducible
$H_i\in\Z[T_1,\ldots,T_n]$ and consider

\begin{equation}\label{xi}
X_i:=\{\nu\ge 1 | \partial_{\nu}(H_i(\mu_k([X_1]),\ldots,\mu_k([X_n])))=0\}
,\; 1\le i\le N.
\end{equation}

By corollary \ref{zeroes2} the $X_i\subset\N$ are good and for fixed
$\nu\ge 1$ we have
\[
0\stackrel{(\ref{relation})}{=}\partial_{\nu}(H(\mu_k([X_1]),\ldots,\mu_k([X_n])))=\prod_{i=1}^N \partial_{\nu}(H_i(\mu_k([X_1]),\ldots,\mu_k([X_n])))\mbox{ in }\Q[T].
\]
As $\Q[T]$ is an integral domain this implies that $\cup_{i=1}^N X_i=\N$.\\
I claim that this implies the existence of $1\le i_0\le N$ and $M\ge 1$ 
such that $M\N\subset X_{i_0}$. Indeed, assume that this is not the case and write
\[
X_i=\Sigma_i\cup_{j\in I_i}(j+M_i\N) \; , \; i=1,\ldots, N
\]

for suitable finite sets $\Sigma_i$, integers $M_i\ge 1$ and subsets
$I_i\subset\{0,\ldots, M_i-1\}$. Then $0\not\in I_i$ for $i=1,\ldots, N$.
Consider a positive multiple $x$ of the product $M_1\dots M_N$ so large that $x\not\in\cup\Sigma_i$.
As $\N=\cup X_i$ we find $i_0$ and $j\in I_{i_0}$ such that $x\in j+M_{i_0}\N$
 hence $x\equiv j\;(\mbox{ mod }M_{i_0})$ and $x\equiv 0\;(\mbox{ mod }M_{i_0})$, a contradiction
because $j\in I_{i_0}$ hence $j\not\equiv 0\;(\mbox{ mod }M_{i_0})$,
so there do indeed exist $M\ge 1$ and $i_0$ such that $M\N\subset X_{i_0}$.
We will see that $G:=H_{i_0}$ and $M$ satisfy the requirements of
theorem \ref{irreducibleeq}.
By construction, $G$ is irreducible and it remains to be checked that
\[
G(\mu_{k_M}([X_1\times_k k_M]),\ldots, \mu_{k_M}([X_n\times_k k_M]))=0.
\]
Equivalently, as $\cap_{\nu\ge 1}\mbox{ker}(\partial_{\nu}:\Rc(\Q)[T]\rightarrow\Q[T])=0$ we need to show that for all $\nu\ge 1$

\begin{equation}\label{compwise}
\partial_{\nu}(G(\mu_{k_M}([X_1\times_k k_M]),\ldots, \mu_{k_M}([X_n\times_k k_M])))=0.
\end{equation}

Using (\ref{33}) and $\partial_{\nu}F_M=\partial_{\nu M}$ \cite{DG}, V, \S 5, 2.6 we obtain
\begin{eqnarray*}
(\ref{compwise}) & = & G(\partial_{\nu}(\mu_{k_M}([X_1\times_k k_M])),\ldots, \partial_{\nu}(\mu_{k_M}([X_n\times_k k_M]))) \\
 & = & G(\partial_{\nu}(F_M(\mu_k([X_1]))),\ldots, \partial_{\nu}(F_M(\mu_k([X_n])))) \\
 & = & \partial_{\nu M}(G(\mu_k([X_1]),\ldots, \mu_k([X_n]))). 
\end{eqnarray*}

However, $\nu M\in M\N\subset X_{i_0}$, so this vanishes by construction
 as $G=H_{i_0}$, see (\ref{xi}).

\end{proofof}

\subsection{Zero dimensional varieties}\label{zerodim}

Let $k$ be a finite field. Before applying the results of \ref{irreducible} we include in this
subsection a digression on the subring of $K_0(\Var_k)$ generated
by zero dimensional varieties.\\
Generally speaking, the difficulty of the problem of deciding whether the classes of given
varieties $[X_1],\ldots ,[X_n]$ are algebraically independent in $K_0(\Var_k)$
grows quite rapidly with $n$. For $n=1$ a complete answer can be given.

\begin{theorem}\label{zerodimthm}

Let $k$ be a finite field and let $X/k$ be separated and of finite type. Then the following are equivalent:
\begin{eqnarray*}
1) & [X]\in K_0(\Var_k) \mbox{ is algebraically {\em dependent}.}\\
2) & [X] \mbox{ is integral over } \Z.\\
3) & \ddim(X)=0.
\end{eqnarray*}
\end{theorem}

Condition 1) means that $[X]$ satisfies some polynomial with coefficients
in $\Z$ but this polynomial does not have to be monic. So, a priori, condition
2) is stronger than condition 1).
Let $\tilde{\Z}\subset K_0(\Var_k)$ denote the integral closure of $\Z$
inside $K_0(\Var_k)$ and let
$S\subset K_0(\Var_k)$ be the subring generated by the classes of zero dimensional varieties.
By theorem \ref{zerodimthm} we have $S\subset\tilde{\Z}$ and even more, for all $X/k$,
if $[X]\in\tilde{\Z}$ then $[X]\in S$. Note that this does not imply 
$\tilde{\Z}\subset S$ because an element of $K_0(\Var_k)$ is in general only a formal difference 
of varieties.\\
In fact, I do not know whether $S\subset\tilde{\Z}$ is an equality.\\
One can give a simple presentation of the ring $S$, see theorem \ref{structureofS} below.

\begin{proofof} theorem \ref{zerodimthm}. We start by showing that 3) implies 2).
Let $X/k$ be zero dimensional, i.e. $X=\Spec(A)$ for some Artinian $k$-algebra
$A$. As $[X]=[X^{red}]$ we can assume that $A$ is reduced, i.e. 
$A=\oplus L_i$ for $L_i|k$ finite field extensions. Then $[X]=\sum 
[\Spec(L_i)]$ and we are reduced to showing that for a finite extension field  $k\subset L$ we have $[\Spec(L)]\in \tilde{\Z}$. But, if $d=[L:k]$, then 

\[
 [\Spec(L)]^2-d\; [ \Spec(L) ]  =  [ \Spec(L\otimes_k L) ] -d\;  [ \Spec(L) ] = 
\]
\[
 [ \Spec(L^{\oplus d}) ] - d\;  [ \Spec(L) ]  =  0.
\]

We used that $L|k$ is Galois to have $L\otimes_k L\simeq L^{\oplus d}$
as $k$-algebras.\\
That 2) implies 1) is trivial and we now show that 1) implies 3).
Let $X/k$ be separated and of finite type and such that $[X]$ is algebraically dependent.
We need to show that dim($X$)=0. As $[X]=[X^{red}]$ and dim($X$)=dim($X^{red}$)
 we assume that $X$ is reduced. As $k$ is perfect this implies
that the smooth locus of $X/k$ is (open and) dense. This will be used later
 but first we exploit the assumption that $[X]$ is algebraically dependent.
There is $0\neq F\in\Z[T_1]$ with $F([X])=0$. Let us consider
$f:=\mu_k([X])\in\Rc(\Q)[T]$. For any $\nu\ge 1$ we have
\[
F(\partial_{\nu}(f))= F(\partial_{\nu}(\mu_k([X])))=(\partial_{\nu}\circ\mu_k)(F([X]))=0.
\]
As $\partial_{\nu}(f)$ lies in $\Q[T]$ and $F\neq 0$ this implies that 
$\partial_{\nu}(f)$ is a constant, i.e. $\frac{d}{dT}(\partial_{\nu}(f))=0$.
Then also
\[
\partial_{\nu}(\frac{df}{dT})=\frac{d}{dT}(\partial_{\nu}(f))=0\mbox{ for all }\nu\ge 1
\]
and using $\cap_{\nu\ge 1}($ker$(\partial_{\nu}:\Rc(\Q)[T]\rightarrow
\Q[T]))=0$ we get $\frac{df}{dT}=0$, i.e. $f$ itself
is a constant, $f\in\Rc(\Q)$.
Now assume by contradiction that $d:=\ddim(X)\ge 1$. 
The coefficient of $T^{2d}$ in $f$ is, by definition of $\mu_k$, 
the image under $\chi_q: K_0(\Rep_{G_k}\Q_l)\hookrightarrow\Lambda(\Q_l)$ of the
virtual Galois representation $\sum_j (-1)^j [\Gr_{2d}^WH^j_c(\oX)]$
and this virtual representation is zero because $f$ is a constant and $2d>0$, i.e.
\begin{equation} \label{virtrep}
\sum_j (-1)^j [\Gr_{2d}^WH^j_c(\oX)]=0.
\end{equation}
From \cite{Weil2}, 3.3.4 we have
\begin{eqnarray}
\Gr_{2d}^WH_c^j(\oX)=0 \mbox{ for } j\neq 2d \mbox{ and} \label{van1}\\
\Gr_{\nu}^WH_c^{2d}(\oX)=0 \mbox{ for } \nu\neq 2d. \label{van2}
\end{eqnarray}

From (\ref{virtrep}) and (\ref{van1}) we get $\Gr_{2d}^WH_c^{2d}(\oX)=0$
which by (\ref{van2}) implies

\begin{equation} \label{vanish}
H_c^{2d}(\oX)=0.
\end{equation}

By the introducing remarks concerning the smooth locus of $X/k$ we can choose
$U\subset X$ open, smooth and purely $d$-dimensional such that the dimension
of $X-U$ is strictly less than $d$ ($U$ need not be dense in $X$; it
may be chosen as the union of the smooth loci of the irreducible 
components of $X$ of dimension $d$). The last condition assures that the
natural map $H_c^{2d}(\overline{U},\Q_l)\rightarrow H_c^{2d}(\oX,\Q_l)$
is an isomorphism hence (\ref{vanish}) gives $H_c^{2d}(\overline{U},\Q_l)=0$. 
Poincar\'e duality gives $H^0(\overline{U},\Q_l)=0$ which is a contradiction because dim$_{\Q_l} H^0(\overline{U},\Q_l)$ is the number of connected
components of $\overline{U}$ and $U\neq\emptyset$.

\end{proofof}

We now give a presentation of the subring $S\subset K_0(\Var_k)$ generated by zero dimensional varieties.
From the proof of theorem \ref{zerodimthm}, 3)$\Rightarrow$2) we know
that $S$ is generated by the classes of finite field extensions of $k$.
Putting $x_{\nu}:=[\Spec(k_{\nu})]\in K_0(\Var_k)$ ($\nu\ge 1$) we have
\[
S=\Z[x_{\nu}:\nu\ge 1]\subset K_0(\Var_k).
\]
For $\mu,\nu\ge 1$ we write $d(\mu,\nu)$ (resp. $m(\mu,\nu)$) for 
the greatest common divisor (resp. the least common multiple) of
$\mu$ and $\nu$. 
Computing tensor products of finite extensions of finite fields one checks that $x_{\mu} x_{\nu}=d(\mu,\nu)x_{m(\mu,\nu)}$ and
we get a surjective ring homomorphism
\begin{equation}\label{Sonto}
\phi:\Z[T_{\nu}:\nu\ge 1]/(T_{\nu}T_{\mu}-d(\mu,\nu)T_{m(\mu,\nu)}:\mu,\nu\ge 1) 
\longrightarrow S,\; T_{\nu}\mapsto x_{\nu}.
\end{equation}

\begin{theorem} \label{structureofS}

The above $\phi$ is an isomorphism.\\
In particular, the $x_{\nu}$ ($\nu\ge 2$) are pairwise
distinct zero divisors in $K_0(\Var_k)$.

\end{theorem}

\begin{rem} The existence of zero divisors in $K_0(\Var_k)$ was 
first addressed by B. Poonen in \cite{Poonen} where he showed their
existence for $k$ of characteristic zero. One can show that for any finitely
generated field $k$ the images of the classes of the non-trivial finite Galois extensions of $k$ in
$K_0(\Var_k)[[\A_k^1]^{-1}]$ are an infinite set of zero divisors.
\end{rem}

\begin{proofof} theorem \ref{structureofS}.
The left hand side of (\ref{Sonto}) is a free $\Z$-module
on the set $\{T_{\nu}:\nu\ge 1\}$.
To show that $\phi$ is injective we use that, for all $n\ge 1$, sending
 $[X]\mapsto |X(k_n)|$ defines a motivic measure $\psi_n:K_0(\Var_k)\rightarrow\Z$
and  that, writing $X_{\nu}:=\Spec(k_{\nu})$, we have

\begin{eqnarray*}
|X_{\nu}(k_n)|=\nu\mbox{ , if }\nu\mbox{ divides }n\\
|X_{\nu}(k_n)|=0  \mbox{ , otherwise,}
\end{eqnarray*}

for all $\nu,n\ge 1$. Now assume that $\phi(\sum_{\nu\ge 1}a_{\nu}T_{\nu})=
\sum_{\nu\ge 1}a_{\nu}[X_{\nu}]=0$
for some $a_{\nu}\in\Z$, almost all zero. We need to show that $a_{\nu}=0$
for all $\nu$. For any $n\ge 1$ we have

\[ \begin{array}{ll}
 (*)_n & 0=\psi_n(\sum_{\nu\ge 1}a_{\nu}[X_{\nu}])=\sum_{\nu\ge 1}a_{\nu}|X_{\nu}(k_n)|=\sum_{\nu|n}a_{\nu}\nu.
\end{array} \]

To show that all $a_{\nu}$ are zero we proceed by induction on the number of prime divisors, counting multiplicity, of $\nu$.\\
For $\nu=1$ we have $(*)_1:a_1=0$.\\
For $\nu\neq 1$ we have from $(*)_{\nu}$ that $\sum_{d|\nu}a_d d=0$. 
By induction hypothesis we have $a_d=0$ for all proper divisors $d$ of
 $\nu$, hence $a_{\nu}\nu=0$ and $a_{\nu}=0$.\\
This shows that $\phi$ is injective, hence an isomorphism. To prove the last 
claim of theorem
\ref{structureofS} we use that,
since $\phi$ is an isomorphism, $S$ is $\Z$-free on the set $\{x_{\nu}:\nu\ge 1\}$ and in
particular that the $x_{\nu}$ are all distinct. Finally, for $\nu\ge 2$,
$0=x_{\nu}^2-\nu x_{\nu}=x_{\nu}(x_{\nu}-\nu)$.
As $x_{\nu}-\nu=x_{\nu}-\nu x_1\neq 0$ (note that $x_1=1$) we see that $x_\nu$ 
is a zero divisor.

\end{proofof}

\begin{rem} \label{example1} 
The above considerations suggest  a simple example illustrating theorem
\ref{irreducibleeq}.
We take $n:=1$, $X_1:=\Spec(k_2)$ in this result and compute 
$\mu_k([X_1])\in\Rc(\Q)[T]$.
We have $H_c^i(\oX_1)=0$ for $i\neq 0$ and $H_c^0(\oX_1)$ is
two dimensional with geometric Frobenius $F_k$ acting by exchanging the
vectors of a (fixed) basis. This gives for $\nu\ge 1$
\begin{equation}\label{trace}
c_{\nu}:= \tr(F_k^{\nu}|H_c^0(\oX_1))=  2 \mbox{ for }\nu\mbox{ even} \\
\end{equation}
and $c_{\nu}=0$ for $\nu$ odd. For arbitrary $G\in\Z[T_1]$ we have
$G(\mu_k([X_1]))=0$ if and only if $G(\partial_{\nu}(\mu_k([X_1])))=0$ 
for all $\nu\ge 1$ and from the above computation of the cohomology
of $X_1$ and proposition \ref{motmeasurefiniteproperties}, i) we get that $\partial_{\nu}(\mu_k([X_1]))\in\Q[T]$ is in fact 
a constant equal to $c_{\nu}$. So $G(\mu_k([X_1]))=0$ is equivalent to
$G(c_{\nu})=0$ for all $\nu\ge 1$. From (\ref{trace}) we see that the 
assertion of theorem \ref{irreducibleeq} holds with $M=2$ and $G=T_1-2$.\\
However, there is no {\em irreducible} polynomial $G\in\Z[T_1]$
satisfying $G(\mu_k([X_1]))=0$ because such a $G$ must, by (\ref{trace}), 
satisfy $G(0)=G(2)=0$ and hence must be divisible by $T(T-2)$.\\ 
This shows that in theorem \ref{irreducibleeq} one has to 
allow for a finite extension of the base field in order to find an
irreducible relation.

\end{rem}

\subsection{Two curves}\label{twocurves}

Let $k$ be a finite field. In this subsection by a curve over $k$ we shall always mean a
proper, smooth and geometrically connected curve. We are going to 
treat in detail the following special case of the problem of algebraic independence 
in $K_0(\Var_k)$.
Let $X$ and $Y$ be curves over $k$. When are $[X]$ and $[Y]$ algebraically independent in 
$K_0(\Var_k)$, i.e. when is the subring $\Z[[X],[Y]]\subset K_0(\Var_k)$
a polynomial ring in the variables $[X]$ and $[Y]$ ?
We will see that this is the case for a generic pair ($X$,$Y$).
First we define special curves which will turn 
out to be those for which our methods can
 shed no light on the above question.
A curve $X/k$ is {\em special} if
i) the degree of $k$ over its prime field is even: $k=\F_q$, $q=p^{2n}$ ($n\ge1$) and
ii) all eigenvalues of the geometric Frobenius attached to $k$
acting on $H^1_c(\oX)$ are equal to $+q^{1/2}$.
Note that by i) $+q^{1/2}$ is a rational integer so ii) is meaningful. 
 Being special is not a property of the curve alone but depends on the
base field, too. The projective line is special over any $k$ satisfying condition i).
An elliptic curve $X/k$ becomes special after a finite extension of 
the base field if and only if it is super-singular. If $X/k$ is 
special then the Jacobian of $X$ is $k$-isogeneous to the power of a
super-singular elliptic curve and we see that most curves of genus at least one
are non-special.
The main result of this subsection is the following.

\begin{theorem}\label{indepoftwocurves}

Let $X_1$ and $X_2$ be proper, smooth and geometrically
connected curves over the finite field $k$ such that
$[X_1]$ and $[X_2]$ are algebraically {\em dependent} in $K_0(\Var_k)$. Then there
is a finite extension $k\subset K$ such that at least one of 
the following holds:\\
1) Both $X_1\times_k K/K$ and $X_2\times_k K/K$ are special.\\
2) The Jacobians of $X_1$ and $X_2$ are $K$-isogeneous.

\end{theorem}

\begin{cor}\label{P1andE}
If $E/k$ is an ordinary elliptic curve then $[\Pe^1]$ and $[E]$
are algebraically independent in $K_0(\Var_k)$.
\end{cor}

This is clear using theorem \ref{indepoftwocurves} and the discussion
preceding it. Theorem \ref{indepoftwocurves} can be applied to the following problem:
It is difficult to exhibit relations in $K_0(\Var_k)$. One such
relation, which follows directly from the definition, is that if
$X\longrightarrow Y$ is a {\em Zariski}-locally trivial fibre bundle
with fibre $G$ then $[X]=[G][Y]$ in $K_0(\Var_k)$. One might hope for a
similar relation for more general $X\longrightarrow Y$, e.g. \'etale 
covers. However:

\begin{cor} Let $X/k$ be a curve of genus at least two which does not
become special after any finite extension of $k$. If $Y\longrightarrow X$
is a non-trivial finite \'etale cover then $[X]$ and $[Y]$ are
algebraically independent in $K_0(\Var_k)$.
\end{cor}

\begin{proof} By the assumption made on $X$, condition 1) of theorem
\ref{indepoftwocurves} is excluded. The Hurwitz formula implies that the 
genus of $Y$ is strictly larger than the genus of $X$, excluding condition
2) of theorem \ref{indepoftwocurves}. Hence theorem \ref{indepoftwocurves}
implies the algebraic independence of $[X]$ and $[Y]$.
\end{proof}

\begin{rem} Rephrasing theorem \ref{indepoftwocurves} we have: If for
two given curves $X_1,X_2/k$ none of 1) or 2) from
 theorem \ref{indepoftwocurves} holds 
after any finite extension of $k$ then $[X_1]$ and $[X_2]$ are algebraically independent.
We will in fact show that in this situation $\mu_k([X_1])$ and $\mu_k([X_2])$
are algebraically independent in $\Rc(\Q)[T]$. We have seen in corollary \ref{H1suffices}
that for the algebraic independence of
the $\mu_k([X_i])$ it is sufficient that
$[H^1_c(\oX_1)]$ and $[H^1_c(\oX_2)]$ are algebraically independent in $K_0(\Rep_{G_k}\Q_l)$.
The following example shows that theorem \ref{indepoftwocurves} is 
sharper than just the algebraic independence of $H^1$'s.
Let $X_2/k$ be a curve such that its Jacobian is $k$-isogeneous
to the square of an ordinary elliptic curve $X_1/k$ (such a curve exists for
suitable $k$ by \cite{BDS}, theorem 2).
Then $H^1_c(\oX_2)\simeq H^1_c(\oX_1)^{\oplus 2}$ as $G_k$-modules
 and so $[H^1_c(\oX_1)]$ and $[H^1_c(\oX_2)]$ are algebraically dependent in $K_0(\Rep_{G_k}\Q_l)$.
But theorem \ref{indepoftwocurves} still shows that the $\mu_k([X_i])$
 ($i=1,2$) (and hence the $[X_i]$ themselves) are algebraically  
independent.
The curve $X_1\times_k K/K$ is not special for any
finite extension $k\subset K$ as explained before theorem
\ref{indepoftwocurves} and the Jacobians of $X_1$ 
and $X_2$ are not isogeneous because they are of
different dimensions.\\
Another consequence of theorem \ref{indepoftwocurves} is that two curves
of different genera have algebraically independent classes unless they both
become special after some finite extension of the base field. As the referee
points out, over $k=\C$ one can easily show that two curves of different
genera have algebraically independent classes using a motivic measure
based on the Hodge polynomial, see \cite{LL2}, 3.4.
\end{rem}

\begin{proofof} theorem \ref{indepoftwocurves}.
Since we are allowed to make a finite extension of $k$ we can assume
by theorem \ref{irreducibleeq} that there is an irreducible $G\in\Z[T_1,T_2]$
with $G(\mu_k([X_1]),\mu_k([X_2]))=0$. Writing
\[
f_{i,\nu}:=\partial_{\nu}(\mu_k([X_i]))\in\Q[T]\; ; \; i=1,2;\nu\ge 1
\]

we have

\begin{equation}\label{Gdoes}
G(f_{1,\nu},f_{2,\nu})=0\mbox{ for all }\nu\ge 1.
\end{equation}

Introducing, for $\nu\ge 1$, 

\begin{eqnarray*}
a_{\nu}:=-\tr(F_k^{\nu}|H^1_c(\oX_1))\\
b_{\nu}:=-\tr(F_k^{\nu}|H^1_c(\oX_2)),
\end{eqnarray*}

where $F_k$ is the geometric Frobenius of $k=\F_q$, say, we have 

\begin{eqnarray*}
f_{1,\nu}=1+a_{\nu}T+q^{\nu}T^2\\
f_{2,\nu}=1+b_{\nu}T+q^{\nu}T^2
\end{eqnarray*}

because $H^0_c(\oX_i)=\Q_l$ and $H^2_c(\oX_i)=\Q_l(-1)$
and using proposition \ref{motmeasurefiniteproperties},i). We consider

\begin{equation}\label{candidate}
F(T_1,T_2):=q(T_2-T_1)^2-(a_1-b_1)(a_1T_2-b_1T_1)+(a_1-b_1)^2.
\end{equation}

The discriminant of $F$, considered as a polynomial in $T_2$, equals

\begin{equation}\label{disc}
4q(a_1-b_1)^2T_1+(a_1-b_1)^2(a_1^2-4q).
\end{equation}

The choice of $F$ is made such that we have:

\begin{eqnarray}\label{candworks}
F(f_{1,1},f_{2,1}) =  F(1+a_1T+qT^2, 1+b_1T+qT^2) & &\\
=q(b_1-a_1)^2T^2-(a_1-b_1)(a_1+a_1b_1T+a_1qT^2- & & \nonumber \\
b_1-b_1a_1T-b_1qT^2)+(a_1-b_1)^2 & & \nonumber \\
=q(b_1-a_1)^2T^2-(a_1-b_1)^2(1+qT^2)+(a_1-b_1)^2 & & \nonumber\\
=(a_1-b_1)^2(qT^2-1-qT^2)+(a_1-b_1)^2=0. & & \nonumber
\end{eqnarray}

To continue, we make the following additional assumption:\\
1) $a_1\neq b_1$:
The prime ideal

\[
P_1:=\kker(\Q[T_1,T_2]\longrightarrow\Q[T],T_i\mapsto f_{i,1})
\]

is proper and non-zero hence of height 1 in $\Q[T_1,T_2]$. By (\ref{Gdoes})
$0\neq G\in P_1$ and $G$ is irreducible. Since $\Q[T_1,T_2]$ is factorial we have
$P_1=(G)$, the principal ideal generated by $G$. 
As $a_1\neq b_1$ the discriminant (\ref{disc}) of $F$
 is not a square in $\Q[T_1]$, so $F$ is irreducible in $\Q[T_1,T_2]$. By (\ref{candworks})
we have $F\in P_1$ and the same reasoning as for $G$ gives $P_1=(F)=(G)$.
So there is some $\alpha\in \Q^*$ with $F=\alpha G$
from which we get
$F(f_{1,\nu},f_{2,\nu})=\alpha G(f_{1,\nu},f_{2,\nu})\stackrel{(\ref{Gdoes})}{=}0$ for all $\nu\ge 1$ which we make more explicit

\begin{eqnarray*}
0=F(f_{1,\nu},f_{2,\nu})=F(1+a_{\nu}T+q^{\nu}T^2,1+b_{\nu}T+q^{\nu}T^2) & & \\
\stackrel{(\ref{candidate})}{=} q((b_{\nu}-a_{\nu})T)^2-(a_1-b_1) & & \\
(a_1+a_1b_{\nu}T+a_1q^{\nu}T^2-b_1-b_1a_{\nu}T-b_1q^{\nu}T^2)+(a_1-b_1)^2 & &\\
=T^2(q(b_{\nu}-a_{\nu})^2-(a_1-b_1)(a_1q^{\nu}-b_1q^{\nu}))+T(-(a_1-b_1) & & \\
(a_1b_{\nu}-b_1a_{\nu}))-(a_1-b_1)^2+(a_1-b_1)^2 & &\\
=T^2(q(b_{\nu}-a_{\nu})^2-q^{\nu}(a_1-b_1)^2)+T((a_1-b_1)(b_1a_{\nu}-a_1b_{\nu})). & & 
\end{eqnarray*}

So for any $\nu\ge 1$ (using $a_1-b_1\neq 0$):

\begin{eqnarray}
q(b_{\nu}-a_{\nu})^2 & = & q^{\nu}(a_1-b_1)^2  \label{first} \\
b_1a_{\nu} & = & a_1b_{\nu}. \label{second}
\end{eqnarray}

Now we assume in addition to $a_1\neq b_1$ that we also have\\
1.1) $a_1\neq 0$ and $b_1\neq 0$:
From (\ref{first}) we have

\[
b_{\nu}^2-2a_{\nu}b_{\nu}+a_{\nu}^2=q^{\nu-1}(a_1-b_1)^2
\]

into which we substitute (\ref{second}) to get 

\[
a_{\nu}^2\frac{b_1^2}{a_1^2}-2a_{\nu}^2\frac{b_1}{a_1}+a_{\nu}^2=q^{\nu-1}(a_1-b_1)^2.
\]

Multiplying this by $a_1^2$ gives

\[
a_{\nu}^2(b_1^2-2a_1b_1+a_1^2)=q^{\nu-1}(a_1-b_1)^2a_1^2
\]

and cancelling $(a_1-b_1)^2\neq 0$ we obtain 

\begin{equation}\label{etappe}
a_{\nu}^2=q^{\nu-1}a_1^2\mbox{ for }\nu\ge 1.
\end{equation}

We have $a_{\nu}^2=\tr(F_k^{\nu}|H^1_c(\oX_1)^{\otimes 2})$, hence

\begin{eqnarray}
\ddet(1-F_kt|H^1_c(\oX_1)^{\otimes 2})\label{quadelle} & = & \eexp\left(-\sum_{\nu\ge 1}a_{\nu}^2\frac{t^{\nu}}{\nu}\right)\\
 \stackrel{(\ref{etappe})}{=} \eexp\left(-\sum_{\nu\ge 1}q^{\nu-1}a_1^2\frac{t^{\nu}}{\nu}\right) & = & \eexp\left(q^{-1}a_1^2\left(\sum_{\nu\ge 1}-\frac{(qt)^{\nu}}{\nu}\right)\right)\nonumber\\
& = & (1-qt)^{a_1^2/q}.\nonumber 
\end{eqnarray}

The first equality in (\ref{quadelle}) is \cite{Hartshorne}, Appendix C, lemma 4.1.
We denote by $g(X_i)$ the genus of $X_i$. As dim$_{\Q_l} H^1_c(\oX_i)=2g(X_i)$ we obtain from (\ref{quadelle}) by comparing degrees:

\begin{equation}\label{platt}
a_1^2=4g(X_1)^2q.
\end{equation}

Denoting by $\alpha_1,\ldots, \alpha_{2g(X_1)}$ the eigenvalues of $F_k$ acting on
$H^1_c(\oX_1)$ considered as complex numbers via some isomorphism
$\overline{\Q_l}\stackrel{\sim}{\rightarrow}\C$ we have 
$|\alpha_i|=q^{1/2}$ and thus

\[
|a_1^2|=|-\sum_{i=1}^{2g(X_1)}\alpha_i|^2\leq(2g(X_1)q^{1/2})^2=4g(X_1)^2q.
\]

By (\ref{platt}) this upper bound is in fact an equality 
which forces all $\alpha_i$ to be equal to either $+q^{1/2}$ or
$-q^{1/2}$. After possibly a quadratic extension of $k$ we will
have all $\alpha_i$ equal to $+q^{1/2}$ and then $X_1/k$ will be
special by definition. As the additional assumptions we have made
so far, namely that $a_1\neq b_1$ and $a_1, b_1\neq 0$, are symmetric
in $X_1$ and $X_2$ the same conclusion holds for $X_2/k$ and we are
in case 1) of theorem \ref{indepoftwocurves}.
Now we assume:\\
1.2) $a_1=0$:\\
We then have, as we still assume that $a_1\neq b_1$, that $b_1\neq 0$ and
(\ref{second}) implies that $a_{\nu}=0$ for all $\nu\ge 1$ hence
$H^1_c(\oX_1)=0$ and $X_1\simeq\Pe^1$. Note that there are
no non-trivial $\overline{k}/k$-forms of $\Pe^1$ because the Brauer group
of a finite field is trivial. From (\ref{first}) we then get
$b_{\nu}^2=q^{\nu-1}b_1^2$
and the same argument as in case 1.1) (c.f. (\ref{etappe})) shows that $X_2$ becomes special after 
at most a quadratic extension. So we are again in case 1) of theorem 
\ref{indepoftwocurves}.\\
1.3) $b_1=0$: this case is symmetric to the case 1.2) above.\\
There is one case left to be considered:\\
2) $a_1= b_1$:\\
We now have $f_{1,1}=f_{2,1}$ and $F:=T_1-T_2$ satisfies
$F(f_{1,1},f_{2,1})=0$. By the same arguments as in case 1)
we get $F(f_{1,\nu},f_{2,\nu})=0$ for all $\nu\ge 1$, i.e.

\[
0=F(1+a_{\nu}T+q^{\nu}T^2,1+b_{\nu}T+q^{\nu}T^2)=(a_{\nu}-b_{\nu})T,
\]

so $a_{\nu}=b_{\nu}$ for all $\nu\ge 1$ from which we get

\[
\ddet(1-F_kt|H^1_c(\oX_1))=\ddet(1-F_kt|H^1_c(\oX_2))
\]

which is equivalent to the Jacobians of $X_1$ and $X_2$ being
$k$-isogeneous \cite{Tate}, theorem 1, c1)$\Leftrightarrow$c2) and hence we are in case 2) of theorem
\ref{indepoftwocurves}.

\end{proofof}

\begin{rem}\label{example2}

We conserve the notations of the above proof.
The simplest case in which theorem \ref{indepoftwocurves}
does not guarantee the algebraic independence of the classes of two curves
in $K_0(\Var_k)$ is if $X_1=\Pe^1$
and $X_2=E$ is a super-singular elliptic curve. Let $k$ be so large that
$X_2/k$ is special and put $q:=|k|$. Then 
$F:=q(T_2-T_1)^2-4qT_1+4q$ satisfies

\begin{equation}\label{nase}
F(\mu_k([X_1]),\mu_k([X_2]))=0.
\end{equation}

In fact, we have 

\begin{eqnarray*}
f_{1,\nu}=1+q^{\nu}T^2 & = & \partial_{\nu}(\mu_k([X_1])) \\
f_{2,\nu}=1-2q^{\nu/2}T+q^{\nu}T^2 & = & \partial_{\nu}(\mu_k([X_2]))
\end{eqnarray*}

and compute

\begin{eqnarray*}
F(\partial_{\nu}(\mu_k([X_1])),\partial_{\nu}(\mu_k([X_2]))) & = & q(-2q^{\nu/2}T)^2-4q(1+q^{\nu}T^2)+4q \\
=4q^{\nu+1}T^2-4q^{\nu+1}T^2-4q+4q & = & 0,
\end{eqnarray*}

for all $\nu\ge 1$, hence (\ref{nase}) follows.
On the other hand

\[
F([X_1],[X_2])=[Z^+]-[Z^-]\mbox{ in }K_0(\Var_k)
\]

where $Z^+$ and $Z^-$ are the (non-connected) proper, smooth varieties

\begin{eqnarray*}
Z^+=(E^2)^{\sqcup q}\sqcup((\Pe^1)^2)^{\sqcup q}\sqcup \Spec(k)^{\sqcup 4q}\\
Z^-=(E\times\Pe^1)^{\sqcup 2q}\sqcup(\Pe^1)^{\sqcup 4q}.
\end{eqnarray*}

The sign $\sqcup$ means disjoint union of schemes.
I cannot decide whether $[Z^+]$ and $[Z^-]$ are equal
(i.e. whether $F([X_1],[X_2])=0$) but I
suspect not and even more, that $[\Pe^1]$ and $[E]$ should be 
algebraically independent.
Note at least that $Z^+\not\simeq Z^-$, e.g. $Z^+$ has isolated points
but $Z^-$ does not.\\
This illustrates that our methods can only give {\em sufficient} 
conditions for varieties to be algebraically independent in the Grothendieck ring
of varieties but we have no idea how to produce interesting relations in
$K_0(\Var_k)$.
\end{rem}

\hspace*{\fill} \\

\end{document}